\documentclass[a4paper,10pt,reqno]{amsart}
\usepackage[margin=2.1cm]{geometry}
\usepackage[utf8]{inputenc}
\usepackage{hyperref}
\usepackage{amsmath}
\usepackage{amsfonts}
\usepackage{amssymb}
\usepackage{amsthm}
\usepackage{latexsym}
\usepackage{xcolor}
\usepackage{tikz}
\usepackage{mathtools}
\usepackage{faktor}
\usepackage{orcidlink}

\providecommand\given{}
\newcommand\SetSymbol[1][]{%
	\nonscript\:#1 :
	\allowbreak
	\nonscript\:
	\mathopen{}}
\DeclarePairedDelimiterX\Set[1]\{\}{%
	\renewcommand\given{\SetSymbol%[\delimsize]
	}
	#1
}

\DeclarePairedDelimiter{\abs}{\lvert}{\rvert}
\DeclarePairedDelimiterXPP{\norm}[2]{}{\lVert}{\rVert}{_{#2}}{#1}
\DeclarePairedDelimiterXPP{\Altnorm}[3]{}{\lVert}{\rVert}{_{#2}^{#3}}{#1}
\DeclarePairedDelimiterXPP{\fNorm}[2]{}{\lvert}{\rvert}{_{#2}}{#1}
\definecolor{darkgreen}{rgb}{0.0, 0.2, 0.13}
\definecolor{darkolivegreen}{rgb}{0.33, 0.42, 0.18}
\definecolor{chamoisee}{rgb}{0.63, 0.47, 0.35}
\definecolor{cerulean}{rgb}{0.0, 0.48, 0.65}
\definecolor{coolgrey}{rgb}{0.55, 0.57, 0.67}

\DeclareMathOperator{\WF}{WF}
\DeclareMathOperator{\real}{Re}
\DeclareMathOperator{\imag}{Im}
\DeclareMathOperator{\supp}{supp}
\DeclareMathOperator{\singsupp}{sing\, supp}
\DeclareMathOperator{\Char}{Char}

\DeclareMathOperator{\Id}{Id}

\DeclareMathOperator{\op}{a-Op}
\DeclareMathOperator{\OPA}{Op}

\newcommand{\fa}{\;\,\forall\,}
\newcommand{\ex}{\;\,\exists\,}

\newcommand{\Cinfty}[1]{\mathcal{C}^\infty\left(#1\right)}

\newcommand{\COT}{T^\ast\Omega\!\setminus\!\{0\}}

\newcommand{\N}{\mathbb{N}}
\newcommand{\R}{\mathbb{R}}

\newcommand{\C}{\mathbb{C}}
\newcommand{\Z}{\mathbb{Z}}

\newcommand{\E}{\mathcal{E}}
\newcommand{\D}{\mathcal{D}}
\newcommand{\CC}{\mathcal{C}}
\newcommand{\An}{\mathcal{C}^\omega}
\newcommand{\G}{\mathcal{G}}
\newcommand{\fG}{\mathfrak{G}}

\newcommand{\bG}{\mathbf{G}}
\newcommand{\bM}{\mathbf{M}}

\newcommand{\bL}{\mathbf{L}}
\newcommand{\bN}{\mathbf{N}}
\newcommand{\fN}{\mathfrak{N}}

\newcommand{\M}{\mathcal{M}}
\newcommand{\eps}{\varepsilon}

\newcommand{\cK}{\mathcal{K}}
\newcommand{\cR}{\mathcal{R}}
\newcommand{\cS}{\mathcal{S}}
\newcommand{\cU}{\mathcal{U}}

\newcommand{\fM}{\mathfrak{M}}

\newcommand{\crb}{\mathcal{V}}

\newcommand{\alp}{{\lvert\alpha\rvert}}
\newcommand{\bet}{{\lvert\beta\rvert}}
\newcommand{\xit}{{\lvert\xi\rvert}}
\newcommand{\etat}{\lvert\eta\rvert}

\DeclareMathOperator{\essupp}{essupp}

\newcommand{\Beu}[2]{\mathcal{E}^{( #1 )} ( #2 )}
\newcommand{\Rou}[2]{\mathcal{E}^{\{ #1 \}} ( #2 )}
\newcommand{\DC}[2]{\mathcal{E}^{[ #1 ]} ( #2 )}

\newcommand{\deriv}[2][]{\frac{\partial^{#1}}{\partial {#2}^{#1}}}
\DeclareMathOperator{\dist}{dist}
\DeclareMathOperator{\bv}{bv}

\newcommand{\dprime}{{\prime\prime}}

\theoremstyle{plain}
\newtheorem{Thm}{Theorem}[section]
\newtheorem{Prop}[Thm]{Proposition}

\newtheorem{Cor}[Thm]{Corollary}
\newtheorem*{MThm}{Theorem}

\theoremstyle{definition}
\newtheorem{Def}[Thm]{Definition}

\theoremstyle{remark}
\newtheorem{Rem}[Thm]{Remark}

\numberwithin{equation}{section}
\allowdisplaybreaks[4]


\begin{document}
	\title[Hypoellipticity in ultradifferentiable classes]{Hypoellipticity of analytic differential operators \\in general ultradifferentiable classes}
	\author{Stefan F\"urd\"os\;\orcidlink{0000-0003-2612-5349}}
	\email{stefan.fuerdoes@univie.ac.at}
	\address{Institute of Mathematics, University of Vienna\\
	Oskar-Morgenstern-Platz 1, 1090 Vienna, Austria}
	%\begin{keywords}
	%	ultradifferentiable hypoellipticity \sep ultradifferentiable wavefront set
	%	\sep analytic pseudodifferential operators \sep (non-)quasianalyticitcity
	%	\sep \MSC[2020]{35B65} 
%	\end{keywords}
		\begin{abstract}
		We show that analytic pseudodifferential and Fourier integral operators behave well for ultradifferentiable classes satisfying minimal regularity properties.
		As an application we investigate the ultradifferentiable regularity properties of several examples of analytic differential operators.
		In particular we extend Treves' characterization of the hypoellipticity of analytic operators of principal type to the ultradifferentiable category.
	\end{abstract}
	
	\maketitle

\section{Introduction}
The focus of this paper is the study of the ultradifferentiable
regularity properties of analytic differential operators
with respect to a large family of ultradifferentiable classes.
In recent years there has been an increased interest in the study of the regularity of partial differential equations
with respect to ultradifferentiable classes,
see e.g.~\cite{Albanese2010}, \cite{Albanese2012}, \cite{BraunRodrigues2021}, \cite{Fuerdoes2024A}, \cite{Hoepfner2020} or\cite{Hoepfner2023}.

We take a slightly different approach compared with the predominant point of view of the literature by asking if the regularity properties of a given analytic partial differential operator $P$
with respect to some ultradifferentiable class influence the regularity properties of $P$ regarding a different ultradifferentiable class.
For example, by the definition of hypoellipticity,  if a differential operator
is hypoelliptic with respect to some ultradifferentiable class
then this knowledge does not in general give information whether
the operator is hypoelliptic in terms of a different class.
But several results, e.g.~\cite{Metivier1980}, \cite{MR2126468},
\cite{Bove2013} or \cite{Cordaro2024},
see also \cite{MR1249275}, indicate
that there is a connection between nondegeneracy conditions on the operators and the property that hypoellipticity with respect to one ultradifferentiable class implies hypoellipticity for any larger class.

 In fact, our first main statement, Theorem \ref{WFkernel}, shows how operators acting on distributions
transforms the wavefront set with respect to suitable general ultradifferentiable classes in terms of the ultradifferentiable wavefront set of their distributional kernel by generalizing a result of H\"ormander \cite{MR0294849} in the case of Denjoy-Carleman classes given by weight sequences
to a wide family of ultradifferentiable classes which includes, among others, ultradifferentiable classes
in the sense of Braun--Meise--Taylor  \cite{MR1052587}.
We might see this also as a continuation of the work done in
 \cite{MR4002151} on the study of microlocal analysis in general ultradifferentiable classes.
In particular,  our statement implies that an analytic parametrix is not only a smooth parametrix but also an ultradifferentiable parametrix.
Another direct consequence is that analytic pseudodifferential operators are microlocal with respect to suitable ultradifferentiable classes and analytic Fourier Integral Operators transform the ultradifferentiable wavefront set
as expected, cf.~Section \ref{Sec:MicroDistr}.

We may refer to ultradifferentiable classes, for which our results hold, as semiregular. For a precise definition of the classes we consider, see the next section and the beginning of Section \ref{Sec:WeightMatrices}.
The results mentioned above induce directly that elliptic analytic pseudodifferential operators are hypoelliptic with respect to any semiregular ultradifferentiable class.
This generalizes well-known results on the ultradifferentiable regularity
of elliptic differential operators with analytic coefficients
by \cite{MR0294849} in the case of Denjoy-Carleman classes, for other classes see \cite{Albanese2010} and \cite{MR4002151}.
But obviously we can also use our methods to prove that other differential operators, which possess an analytic parametrix, are actually ultradifferentiable hypoelliptic.
In particular, we are able to extend Treves' characterization
of hypoellipticity of analytic differential operators
of principal type to the ultradifferentiable setting, which
can formulated slightly loosely in the following way:
\begin{MThm}
	Let $P$ be an analytic differential operator of principal type
	and $\mathcal{U}$ be a semiregular ultradifferentiable class. 
	Then $P$ is smooth hypoelliptic if and only if
	$P$ is hypoelliptic with respect to $\mathcal{U}$.
\end{MThm}
For a precise formulation see Theorem \ref{PrincipalThm1}
and Theorem \ref{GeneralTrevesThm}.
Other operators we are looking at are those considered in
\cite{MR2126468} and in \cite{MR0597752}.

Therefore our results, summarized in the next section, show that
operators with very nice nondegeneracy properties have nice
regularity properties for a wide family of ultradifferentiable classes. 
Using more general ultradifferentiable classes 
in the study of regularity of PDOs might give more precise information than the usual approach with Gevrey classes.
In the rest of the introduction we try to give an overview
of the existing literature in order to indicate
how the results of this paper fit into it.
We may start by noticing that analytic hypoellipticity of
operators does not generally imply hypoellipticity in the smooth category:
For example
\begin{equation*}
	P=\left(D_1+ix_1^{2\ell}D_2\right)^2+cD_2,
\end{equation*}
where $\ell\in\N$ and $c\in\C\!\setminus\!\{0\}$,
is analytic hypoelliptic but not hypoelliptic,
cf.~\cite{Okaji1988} and for other examples in this regard see \cite{MR1249275}.
In fact, \cite{Okaji1988} showed that $P$ is
$\G^s$-hypoelliptic for $1\leq s\leq 4\ell/(2\ell-1)$
where $\G^s$ denotes the Gevrey class of degree $s\geq 1$, 
i.e.~the ultradifferentiable class associated to the weight sequence $k!^s$, cf.~Definition \ref{DC-Def}.
In contrast, it is a classical statement that
the heat operator $\partial_t-\Delta_x$ is
smooth and $\G^s$-hypoelliptic if and only if $s\geq 2$.
More generally the characterization of hypoellipticity of differential operators with constant coefficients \cite{Hoermander1955}, see also \cite{Hoermander2005}
and \cite{MR1249275}, implies that this is generally true for hypoelliptic operators $P$ with constant coefficients:
If $P$ is smooth hypoelliptic then there is some $s_0\geq 1$
such that $P$ is $\G^s$-hypoelliptic for $s\geq s_0$, where
$s_0$ depends on the geometry of the symbol $P(\xi)$ of $P$. 
This indicates that there might be a connection between 
the geometry of the symbol of some operator with its hypoelliptic
properties with respect of Gevrey or, more generally, ultradifferentiable classes.

Historically, beginning with the study of differential operators
with constant operators,  the Gevrey classes $\G^s$ have been
 used to interpolate between the algebra of smooth functions $\CC^\infty$ and the space of real-analytic functions $\CC^\omega$, see \cite{MR1249275} and the references therein.
A current example of this approach is the ongoing investigation
of the problem of analytic hypoellipticity of 
sum of squares of analytic real vector fields which satisfy the
H\"ormander condition, see for example \cite{Albano2013},
\cite{Bove2018}, \cite{Bove2024}, \cite{Chinni2023}.
The main reason, why this approach is sensible is that
if $P$ is an analytic sum of squares operator satisfying 
the condition of H\"ormander 
then there is some $s\geq 1$ such that $P$ is $\G^s$-hypoelliptic, cf.~\cite{Derridj1973} and if  $P$ is $\G^s$-hypoelliptic then
$P$ is $\G^t$-hypoelliptic for any $t\geq s\geq 1$.
The last statement is actually a corollary of Metivier's inequality, which in turn characterizes analytic and Gevrey hypoellipticity for operators which are $L^2$-solvable, see
\cite{Metivier1980} and \cite{Bove2013}.

With respect to ultradifferentiable classes
Metivier's inequality has been generalized by Cordaro and the author \cite{Cordaro2024} to Denjoy-Carleman classes $\E^{\{\bM\}}$ given
by admissible weight sequences $\bM$. Here admissible weight sequences
are semiregular (see Definition \ref{Def:Semiregular})
and satisfy the following condition
\begin{equation}\label{mg}\tag{$\star$}
	\ex Q>0:\; M_{k+\ell}\leq Q^{k+\ell+1}M_kM_\ell,\qquad k,\ell\in\Z_+.
\end{equation}
The important point here is that we do not require in \cite{Cordaro2024} that $\Rou{\bM}{\Omega}$ is non-quasianalytic.
Recall that a subalgebra $E$ of $\CC^\infty(\Omega)$ is non-quasianalytic if there are non-trivial ultradifferentiable functions with compact support, i.e.~$E\cap \CC^\infty_0(\Omega)\neq \{0\}$.
The well-known Denjoy-Carleman Theorem (see e.g.~\cite{MR1996773}) states that $\Rou{\bM}{\Omega}$
is non-quasianalytic if and only if
\begin{equation}\label{NQ}\tag{$\mathbf{nq}$}
	\sum_{k=1}^\infty\frac{M_{k-1}}{M_k}<\infty.
\end{equation}
We say that the weight sequence $\bM$ is quasianalytic 
if \eqref{NQ} is not satisfied.

As in the Gevrey setting it is a direct consequence of the ultradifferentiable Metivier inequality that $\{\bM\}$-hypoellipticity, $\bM$ being admissible, of
an $L^2$-solvable differential operator implies $\{\bN\}$-hypoellipticity for any larger class $\Rou{\bN}{\Omega}$, where $\bN$ is admissible, and smooth hypoellipticity.
But we need to point out that admissibility is a rather strict
condition for weight sequences. More precisely, if $\bM$
satisfies \eqref{mg} then there is some $s>1$ such that
$\E^{\{\bM\}}\subseteq\G^s$, cf.~\cite{Matsumoto1987}.
This means that many reasonable Denjoy-Carleman classes
are not admissible:
In particular the results in \cite{Cordaro2024} do not cover
the classes given by the weight sequences
$\bN^q=q^{k^2}$, $q>1$, which are semiregular (cf. Definition \ref{Def:Semiregular}).
Another interesting example of a semiregular weight sequence, which is not admissible, is given in \cite{MR4002151}:
There exists a semiregular weight sequence $\bL$ with the following properties:
\begin{itemize}
	\item $\bL$ is quasianalytic.
	\item $\E^{\{\bM\}}\nsubseteq\G^s$ for all $s>1$.
\end{itemize}

Hence using the ultradifferentiable Metivier inequality we could not infer that an analytic hypoelliptic differential operator of principal type is hypoelliptic with respect to the classes given by either the sequence $\bN$ or $\bL$, which is
true due our main Theorem, cf. Theorem \ref{PrincipalThm1}.
But in fact the setting we work in this paper is still much larger regarding the scope of ultradifferentiable classes we consider:
First, in \cite{Cordaro2024} we only worked with Roumieu classes
given by weight sequences, but here we consider also Beurling
classes (see Definition \ref{DC-Def} for the difference).

Moreover, as  alluded to above, the setting we work with in this paper includes ultradifferentiable classes given by
weight functions in the sense of Braun--Meise--Taylor \cite{MR1052587}.
Gevrey classes can be realized both as Denjoy-Carleman classes
and as Braun-Meise-Taylor classes, but in general the classes in the sense of Braun--Meise--Taylor classes cannot described
as Denjoy-Carleman classes and vice-versa, cf.~\cite{Bonet2007}.
In fact, we work in the setting of ultradifferentiable classes
given by weight matrices, i.e.~families of weight sequences,
introduced by \cite{MR3285413}. These classes include both
Denjoy-Carleman classes and Braun-Meise-Taylor classes and using
weight matrices allows to take an unified approach.
We may remark that the advantage to work in the setting of weight matrices is
particularly apparent in the study of the Problem of Iterates, see \cite{Fuerdoes2022}. Here we note that we can construct
semiregular quasianalytic classes given by weight matrices, which are not
contained in any semiregular Denjoy-Carleman classes, 
cf.~Corollary \ref{TransversalCor}.

The paper is organized in the following way: In Section \ref{Sec:First} the results on the ultradifferentiable hypoellipticity are summarized for easier reference in terms of Denjoy-Carleman and Braun-Meise-Taylor classes.
In Section \ref{Sec:WeightMatrices} we first recall 
the basic facts about the ultradifferentiable wavefront set
associated to classes given by weight matrices from \cite{MR4002151} and continue to prove the generalization of
H\"ormander's theorem on the ultradifferentiable microlocal transform properties of operators in terms of their distributional kernel.

In Section \ref{Sec:MicroDistr} we apply this theorem to
analytic pseudodifferential operators and analytic Fourier integral operators. In particular we prove the microlocal regularity theorem for classical analytic pseudodifferential operators in the category of semiregular ultradifferentiable classes.
Section \ref{Sec:principal} is devoted to the generalization
of Treves' Theorem on hypoelliptic differential operators of
principal type.
Finally, in Section \ref{Sec:Metivier} we prove the ultradifferentiable version of the main result of \cite{MR0597752}.

We have included two appendices, which might be of interest.
In the first appendix we quickly summarize the basic theory on
ultradifferentiable hypoellipticity of differential operators with constant coefficients in terms of classes given by weight matrices. 
 This might be well-known to experts in the field, but 
it seems not to have been written down in the literature beyond the Gevrey case.

The second appendix is basically devoted to show that there
is a quasianalytic semiregular ultradifferentiable class given by a weight matrix (even a weight function) which is \emph{not} contained in any
semiregular Denjoy-Carleman class.

\subsection*{Acknowledgments}\leavevmode

This work was funded in whole or in part by the Austrian Science Fund (FWF) 10.55776/PAT1994924.
The starting point of this article was a question that Paulo Cordaro posed, while the author was a Post-Doc under his supervision at the University of S\~{a}o Paulo in 2023.

%The author would also like to thank Francois Treves and
%Bernhard Lamel for their interest in this work.

Finally the author expresses his gratitude to Gerhard Schindl
for pointing out and allowing to include the proof of Theorem \ref{TransversalThm} in this paper.

\section{Statement of main results}\label{Sec:First}
%\begin{Def}
%	with symbol $p\sim \sum_j p_{d-j}$. The characteristic set of $P$ is
%	\begin{equation*}
%		\Char P=\Set*{(x,\xi)\in\COT\given p_{d}(x,\xi)=0}.
%	\end{equation*}
%\end{Def}

\begin{Def}\label{DC-Def}
	A  sequence $\bM=(M_k)_k$ of positive numbers is called a weight sequence if $M_0=1$,
	$\sqrt[k]{M_k}\rightarrow\infty$ and 
	\begin{equation}\label{LogConvex}
		M_k^2\leq M_{k-1}M_{k+1},\qquad k\in\N.
	\end{equation}
The Roumieu class $\Rou{\bM}{\Omega}$ associated to $\bM$ on an open set $\Omega\subseteq\R^n$
consists of the smooth functions $f\in\CC^\infty(\Omega)$ such that for each compact set $K\subseteq\Omega$
there are constants $C,h>0$ so that the following estimate holds:
\begin{equation}\label{DenjoyDef}
	\sup_{x\in K}\,\abs*{D^\alpha f(x)}\leq Ch^\alp M_\alp,\qquad \fa \alpha\in\Z_+^n.
\end{equation}
On the other hand a smooth function $f\in\CC^\infty(\Omega)$ is an element of the Beurling class
$\Beu{\bM}{\Omega}$ associated to $\bM$ if for all compact sets $K\subseteq\Omega$ and every
$h>0$ there is a constant $C>0$ such that \eqref{DenjoyDef} is satisfied.
\end{Def}

\begin{Def}\label{Def:BMT}
	A weight functions $\omega$ is a non-decreasing continuous function 
	$\omega: [0,\infty)\rightarrow [0,\infty)$ such that $\omega (t)=0$ for $t\in[0,1]$ and
	the following conditions hold:
	\begin{align}\tag{$\omega_1$}\label{omega1}
		\omega(2t)&\,=O(\omega(t)),\qquad t\rightarrow\infty,\\
\tag{$\omega_2$}		\log t&\,=o(\omega(t)),\qquad t\rightarrow \infty,\\
	\tag{$\omega_3$}	\varphi_\omega(t)&:=\omega(e^t)\;\text{is convex}.
	\end{align}
We need also the conjugate function of $\varphi_\omega$:
\begin{equation*}
\varphi^\ast_\omega(t)=\sup_{s\geq 0} (st-\varphi_\omega(s)).
\end{equation*}
The Roumieu class $\Rou{\omega}{\Omega}$ associated to $\omega$ is the space of functions $f\in\CC^\infty(\Omega)$ 
such that for every compact set $K\subseteq\Omega$ there are constants $C,h>0$ so that the following estimate
is satisfied
\begin{equation}\label{BMT-Est}
	\sup_{x\in K}\abs*{D^\alpha f(x)}\leq Ce^{\varphi^\ast_\omega(h\alp)},\quad \alpha\in\Z_+^n.
\end{equation}
A function $f\in\CC^\infty(\Omega)$ is an element of the Beurling class $\Beu{\omega}{\Omega}$
if for each compact set $K\subseteq\Omega$ and every $h>0$ there is a constant $C>0$ such
that \eqref{BMT-Est} holds.
\end{Def}

We will use the notation $[\ast]=\{\ast\},(\ast)$ ($\ast=\bM,\omega$ etc.) if the concerned statements hold
for Roumieu and Beurling classes alike.

\begin{Def}\label{Def:Semiregular}\leavevmode
	\begin{enumerate}
		\item A weight sequence $\bM$ is semiregular if the following conditions hold:
		\begin{gather}
		\label{DC-AnalyticIncl}	\lim_{k\rightarrow\infty}\frac{\sqrt[k]{M_k}}{k}=\infty,\\
			\ex Q>0:\quad M_{k+1}\leq Q^{k+1}M_k, \qquad \fa k\in\Z_+.\label{DC-DerivClosed}
		\end{gather}
		\item We say that a weight function is semiregular if 
	\end{enumerate}

\begin{equation}\label{BMT-AnalyticIncl}
	\omega(t)=o(t),\qquad t\rightarrow\infty.
\end{equation}
\end{Def}

\begin{Rem}
	If $\bM$ (resp.\ $\omega$) is semiregular then $\DC{\bM}{\Omega}$ (resp. $\DC{\omega}{\Omega}$)
	is closed under derivation and contains strictly the class $\An(\Omega)$ of analytic functions:
	In the case of weight sequences condition \eqref{DC-AnalyticIncl} gives instantly that
	$\An(\Omega)=\G^1(\Omega)\subsetneq\DC{\bM}{\Omega}$. 
	Moreover, $\eqref{DC-DerivClosed}$ implies that if $f\in\DC{\bM}{\Omega}$
	then $D^\alpha f\in\DC{\bM}{\Omega}$ for any $\alpha\in\Z_+^n$.
	
	Regarding classes given by weight functions, derivation closedness of $\DC{\omega}{\Omega}$
	follows from the convexity of $\varphi_\omega$ , cf.~\cite{MR1052587}.
	On the other hand, \eqref{BMT-AnalyticIncl} gives that
	$\An(\Omega)\subsetneq\DC{\omega}{\Omega}$, see \cite{MR3285413}.
\end{Rem}

\begin{Def}%\label{Def:SemiregularPseudo}
	Let $P$ be a pseudodifferential operator in $\Omega\subseteq\R^n$.
	\begin{enumerate}
		\item If $\bM$ is a weight sequence then we say that $P$ is $[\bM]$-hypoelliptic
		if 
		\begin{equation*}
			\singsupp_{[\bM]} Pu=\singsupp_{[\bM]} u
		\end{equation*}
	for all $u\in\E^\prime(\Omega)$.
	\item If $\omega$ is a weight function then $P$ is $[\omega]$-hypoelliptic if
	\begin{equation*}
		\singsupp_{[\omega]}Pu=\singsupp_{[\omega]} u
	\end{equation*}
for all $u\in\E^\prime(\Omega)$.
	\end{enumerate}
\end{Def}
Here $\singsupp_{[\ast]}u$ denotes the ultradifferentiable singular support of a distribution $u$
which is defined in the obvious way.
We can now state the main results of this article. 
The first one extends statements of \cite{MR0294849} in the Roumieu case, for the Beurling case
see \cite{MR4002151},
for analytic partial differential operators.
%For the definition of the ultradifferentiable wavefront set of a distribution we refer to
%the next section.
\begin{Thm}\label{StandardElliptRegThm}
	Let $A$ be an elliptic analytic pseudodifferential operator in $\Omega$.
	Then $A$ is $[\bM]$-hypoelliptic for every semiregular
	weight sequence and $[\omega]$-hypoelliptic for all
	semiregular weight functions
	%and $u\in\E^\prime(\Omega)$.
	%Then the following holds:
	%\begin{enumerate}
	%	\item If $\bM$ is a semiregular weight sequence then
	%	\begin{equation*}
	%		\WF_{[\bM]}Au\subseteq \WF_{[\bM]} u\subseteq\WF_{[\bM]} Au\cup\Char A.
	%	\end{equation*}
%	\item If $\omega$ is a weight function such that $\omega(t)=O(t)$ then
	%\begin{equation*}
	%	\WF_{[\omega]}Au\subseteq\WF_{[\omega]}u\subseteq %\WF_{[\omega]}Au\cup \Char A
	%\end{equation*}
	%\end{enumerate}
\end{Thm}
%\begin{Cor}
%	Let $A$ be an elliptic classical analytic pseudodifferential operator in $\Omega$.
%	Then $A$ is $[\bM]$-hypoelliptic for all semiregular weight sequence $\bM$
%	and $[\omega]$-hypoelliptic for any semiregular weight function $\omega$.
%\end{Cor}

The simplest class of non-elliptic operators are the operators of principal type:
\begin{Def}
	A differential operator $P$ of order $d$ with principal symbol $p_d$ on $\Omega$ is of principal type\footnote{We use here the classical definition of \emph{principal type} as used in \cite{MR0296509}, which differs from the modern definition, e.g.~\cite[Definition 23.1.2]{MR4436039}} if 
	\begin{equation*}
		\abs*{p_d(x,\xi)}+\abs*{d_\xi p_d(x,\xi)}\neq 0
	\end{equation*}
for all $(x,\xi)\in\COT$.
\end{Def}

Treves \cite{MR0296509} proved that a differential operator of principal type with analytic coefficients is hypoelliptic if and only if it is analytic hypoelliptic.
We can extend this equivalence to include any semiregular ultradifferentiable class:
\begin{Thm}\label{PrincipalThm1}
	Let $A$ be a differential operator of principal type with analytic coefficients,
	$\bM$ be a semiregular weight sequence and $\omega$ be a semiregular weight function.
	Then the following statements are equivalent:
	\begin{enumerate}
		\item $A$ is smooth hypoelliptic.
		\item $A$ is $\{\bM\}$-hypoelliptic.
		\item $A$ is $(\bM)$-hypoelliptic.
		\item $A$ is $\{\omega\}$-hypoelliptic.
		\item $A$ is $(\omega)$-hypoelliptic.
		\item $A$ is analytic hypoelliptic.
	\end{enumerate}
\end{Thm}
Our methods also allow us to generalize a theorem of Metivier \cite{MR0597752}
on the analytic hypoellipticity of operators with multiple characteristics (for the analogous result regarding smooth hypoellipticity see  \cite{BoutetdeMonvel1976}):
\begin{Thm}\label{MetivierThm}
		Let $P$ be a classical analytic pseudodifferential operator of order $m$ with symbol
	\begin{equation*}
		p(x,\xi)\sim\sum_{j=0}^\infty p_{m-j}(x,\xi)
	\end{equation*}
	and suppose that $\Char P$ is a symplectic analytic submanifold of $T^\ast\Omega\!\setminus\!\{0\}$.
	Assume further that there is some $k\in\N$ such that $p_m$ vanishes exactly of order $k$ on $\Char P$
	and that $p_{m-j}$ vanishes of order $k-2j$ on $\Char P$ for $j\leq k/2$.
	If finally $P$ is $k/2$-subelliptic then $P$ is $[\bM]$-hypoelliptic in $\Omega$ for any
	semiregular weight sequence $\bM$ and $[\omega]$-hypoelliptic for all semiregular
	weight functions $\omega$.
\end{Thm}
Finally our approach also leads directly to an extension of  the main statements in \cite{MR2126468}:
Let $k\in\N$ and $\lambda\in\C$. We consider the linear differential operator
\begin{equation*}
	Q_\lambda=\deriv[2]{t}+t^{2k}\deriv[2]{x}-ik(1+2\lambda)t^{k-1}\deriv{x}.
\end{equation*}

\begin{Thm}\label{Cordaro1}
	Let $k\geq 2$ be an even integer,
	$\bM$ be a semiregular weight sequence, $\omega$ a weight function with $\omega(t)=o(t)$,
	 $\Omega$ be an open set with $\Omega\cap\{t=0\}\neq\emptyset$.
	
	Then the following statements are equivalent:
	\begin{enumerate}
		\item $Q_\lambda$ is hypoelliptic in $\Omega$.
		\item $Q_\lambda$ is $\{\bM\}$-hypoelliptic in $\Omega$.
		\item $Q_\lambda$ is $(\bM)$-hypoelliptic in $\Omega$.
		\item $Q_\lambda$ is $\{\omega\}$-hypoelliptic in $\Omega$.
		\item $Q_\lambda$ is $(\omega)$-hypoelliptic in $\Omega$.
		\item $Q_\lambda$ is analytic hypoelliptic in $\Omega$.
		\item $\lambda\notin \Lambda_k=\Set*{j(1+1/k)+1/(2k)\given j\in\Z}$.
	\end{enumerate}
\end{Thm}

\section{Microlocal analysis in general ultradifferentiable classes}\label{Sec:WeightMatrices}
%We assemble here the basic definition about the ultradifferentiable wavefront set.
%In order to streamline the presentation of the proofs we will work from now in the setting 
%of weight matrices, which generalizes the theory of weight sequences and weight functions.
\begin{Def}
	A weight matrix $\fM$ is a collection of weight sequences such that
	for any pair $\bM,\bN\in\fM$ we have either $M_j\leq N_j$ for all $j\in\Z_+$
	or $N_j\leq M_j$ for all $j\in \Z_+$.
\end{Def}
If $\Omega\subseteq\R^n$ is an open set then a smooth function $f\in\Cinfty{\Omega}$ is
an element of the Roumieu class $\Rou{\fM}{\Omega}$ associated to the weight matrix $\fM$ if
for all compact sets $K\subseteq\Omega$ there are constants $C,h>0$ and some weight sequence
$\bM\in\fM$ such that
\begin{equation}\label{RSestimate}
	\sup_{x\in K}\abs*{D^\alpha f(x)}\leq Ch^\alp M_\alp\qquad \fa\alpha\in \Z_+^n.
\end{equation}
On the other hand $f$ is in the Beurling class $\Beu{\fM}{\Omega}$ associated to $\fM$ 
if for every compact set $K\subseteq\Omega$, every $\bM\in\fM$ and all $h>0$ there is
a constant $C>0$ such that \eqref{RSestimate} is satisfied.

%\begin{Rem}
%	Let $\bM$ be a weight sequence and $h>0$. If $U\subseteq\R^n$ is an open subset
%	then
%	\begin{equation*}
%		\E^{\bM}_h(U)=\Set*{f\in\CC^\infty(U)\given \norm*{f}{U,\bM,h}
%			=\sup_{\substack{x\in U\\\alpha\in\Z_+^n}}\frac{D^\alpha f(x)}{h^\alp M_\alp}<\infty}
%	\end{equation*}
%is a Banach space.
%We recall that according to \cite{MR0320743}
%the canonical embedding $\Lambda^h_{h_1}$, $0<h_1\leq h$,
%\begin{equation*}
%	\E^{\bM}_{h_1}(U)\longrightarrow\E^{\bM}_{h}(U)
%\end{equation*}
%is compact. Then the spaces
%\begin{align*}
%	\Rou{\fM}{\overline{U}}=\Set*{f\in \CC^\infty(U)\given \ex\bM\in \fM\;\ex h>0:\quad f\in\E^{\bM}_h}\\
%	\shortintertext{and}
%	\Beu{\fM}{\overline{U}}=\Set*{f\in\CC^\infty(U)\given \fa \bM\in\fM;\fa h>0:\quad
%	f\in\E^{\bM}_h}
%\end{align*}
%are DFS spaces and Frechet spaces, respectively.
% \end{Rem}
\begin{Def}
	Let $\fM$ be a weight matrix and $P$ be a linear differential operator with
	coefficients in $\DC{\fM}{\Omega}$. 
	Then we say that $P$ is $[\fM]$-hypoelliptic at $x_0\in\Omega$ if
	there is a neighborhood $U_0$ of $x_0$ such that
	\begin{equation*}
		\singsupp_{[\fM]}Pu=\singsupp_{[\fM]} u
	\end{equation*}
for all $u\in\D^\prime(U_0)$.

The operator $P$ is $[\fM]$-hypoelliptic in an open set $U\subseteq\Omega$ if $P$ is $[\fM]$-hypoelliptic at every point $x$ of $U$.
\end{Def}
\begin{Def}
	Let $\fM$ be a weight matrix.
	\begin{enumerate}
		\item We say that $\fM$ is $R$-semiregular if
		\begin{gather}
			\label{AnalIncl}	\fa\bM\in\fM\;\,\lim_{j\rightarrow\infty}\left(\frac{M_j}{j!}\right)^{1/j}=\infty\\
			\shortintertext{and}
\label{R-DerivClosed}			
\fa \bM\in\fM\,\ex\bN\in\fM\,\ex Q>0:\;\, M_{j+1}\leq Q^{j+1}N_j\qquad \fa j\in\Z_+
		\end{gather}
		hold.
		\item The weight matrix $\fM$ is $B$-semiregular if \eqref{AnalIncl} and
		\begin{equation}\label{B-DerivClosed}
			\fa\bN\in\fM\,\ex\bM\in\fM\,\ex Q>0:\;\, M_{j+1}\leq Q^{j+1}N_j\qquad \fa j\in\Z_+
		\end{equation}
		are satisfied.
	\end{enumerate}
\end{Def}	We write [semiregular]$=R$-semiregular,$B$-semiregular.
If $\fM$ is [semiregular] then $\DC{\fM}{\Omega}$ contains $\An(\Omega)$ by \eqref{AnalIncl} and is closed
under derivation.
In fact, $\DC{\fM}{\Omega}$ is closed under real-analytic mappings, cf.~\cite{MR4002151}.

%We close this section to give a criterion, which will help us to prove non-hypoellipticity.
%It is well-known that (smooth, analytic or ultradifferentiable) hypoellipticity implies
%certain a-priori estimates, see e.g.\  \cite{MR4436039} and in a more general context %\cite{MR1916845}.
%We prove here a stronger result which was inspired by \cite[Lemma 3.1]{MR3043156}.

%\section{Microlocal analysis in ultradifferentiable classes}
%\subsection{The ultradifferentiable wavefront set}
We now recall the microlocalization of the concept of ultradifferentiability from \cite{MR4002151}
\begin{Def}
	Let $\fM$ be a weight matrix and $u\in\D^\prime(\Omega)$.
	\begin{enumerate}
		\item $u$ is microlocally ultradifferentiable of class $\{\fM\}$
		at $(x_0,\xi_0)\in T^\ast\Omega\!\setminus\!\{0\}=\Omega\times\R^n\!\setminus\!\{0\}$
		if there exists a neighborhood $V$ of $x_0$, a conic neigborhood $\Gamma$ of $\xi_0$
		and a bounded sequence $(u_j)_j\subseteq\E^\prime(\Omega)$ with $u_j\vert_V=u\vert_V$
		such that for some $\bM\in\fM$ and some $q>0$ we have that
		\begin{equation}\label{WFestimate}
			\sup_{\substack{\xi\in\Gamma\\ j\in\Z_+}}\frac{\xit^j\abs*{\hat{u}_j(\xi)}}{q^j M_j}<\infty.
		\end{equation}
		\item $u$ is microlocally ultradifferentiable of class $(\fM)$ at $(x_0,\xi_0)$ if
		there are a neighborhood $V$ of $x_0$, a conic neighborhood $\Gamma$ of $\xi_0$ 
		and a bounded sequence $(u_j)_j\subseteq\E^\prime(\Omega)$ with $u_j\vert_V=u\vert_V$
		such that for all $\bM\in\fM$ and all $q>0$ the estimate \eqref{WFestimate} holds.
	\end{enumerate}
	We set
	\begin{equation*}
		\WF_{[\fM]}u:=\Set*{(x,\xi)\in T^\ast\Omega\!\setminus\!\{0\}\given
			u\text{ is not microlocally ultradiff.~of class }[\fM]\text{ at }(x,\xi)}.
	\end{equation*}
\end{Def}
The basic properties of $\WF_{[\fM]}u$ (cf.~\cite[Proposition 5.4]{MR4002151}) are
\begin{Prop}\label{WFproperties}
	Let $\fM$ be a weight matrix and $u\in\D^\prime(\Omega)$.
	\begin{enumerate}
		\item $\WF_{[\fM]}u$ is a conic and closed subset of $T^\ast\Omega\setminus\{0\}$.
		\item $\WF_{\{\fM\}}u\subseteq\WF_{(\fM)}u$.
		\item If $\fM$ is [semiregular] then
		\begin{equation*}
			\pi_1\left(\WF_{[\fM]}u\right)=\singsupp_{[\fM]}u.
		\end{equation*} 
		\item If $\fM$ satisfies \eqref{AnalIncl} then 
		$\WF u\subseteq\WF_{[\fM]}u\subseteq\WF_a u$.
		\item If $\fM$ is [semiregular] and $P$ is a linear differential operator with 
		coefficients in $\DC{\fM}{\Omega}$ then $\WF_{[\fM]}Pu\subseteq\WF_{[\fM]}u$.
	\end{enumerate}
\end{Prop}
Following the arguments in \cite[Section 8.6]{MR1996773} we can prove the microlocal elliptic regularity theorem in $\E^{[\fM]}$ for analytic differential operators:
\begin{Thm}\label{Microelliptic}
	Let $P$ a differential operator with analytic coefficients in $\Omega$ and $\fM$ be
	a [semiregular] weight sequence. Then
	\begin{equation*}
		\WF_{[\fM]} u\subseteq\WF_{[\fM]}Pu\cup\Char P
	\end{equation*}	
	for all $u\in\D^\prime(\Omega)$.
\end{Thm}
For a more general version for operators with ultradifferentiable coefficients see \cite[Section 7]{MR4002151}.
As pointed out above we restrict ourselves to analytic differential operators. Nevertheless,
for our purposes we still need to extend the microlocal theory presented in \cite{MR4002151}
by generalizing some statements from \cite{MR1996773} which have not been considered in \cite{MR4002151}.

We recall that the tensor product of two distributions $u\in\D^\prime(U)$ and $v\in\D^\prime(V)$
($U\subseteq\R^m$, $V\subseteq\R^n$ open), denoted by $u\otimes v$, is defined by
\begin{equation*}
	\langle u\otimes v,\varphi\otimes\psi\rangle=\langle u,\varphi\rangle\cdot\langle v,\psi\rangle
\end{equation*}
where $\varphi\in\CC^\infty_0(U)$ and $\psi\in\CC^\infty_0(V)$, cf.~\cite[Chapter 5]{MR1996773}.
\begin{Thm}\label{WFTensor}
	Let $U\subseteq\R^m$ and $V\subseteq\R^n$ be open sets and $\fM$ be a 
	[semiregular] weight matrix.
	Then
	\begin{equation*}
		\begin{split}
		\WF_{[\fM]}(u\otimes v)&\subseteq \left(\WF_{[\fM]}u\times\WF_{[\fM]}v\right)
		\cup\left(\left(\supp u\times\{0\}\right)\times\WF_{[\fM]}v\right)\\
		&\qquad\cup\left(\WF_{[\fM]}u\times\left(\{0\}\times\supp v\right)\right).
		\end{split}
	\end{equation*}
\end{Thm}
\begin{proof}
	Since $\singsupp_{[\fM]} u\subseteq\supp u$ it is enough to show that if
	\begin{equation*}
		\begin{split}
			(x_0,\xi_0,y_0,\eta_0)&\in \left(\left(U\!\setminus\!(\supp u)\times\R^m\right)
			\times V\times\R^n\right)
			\cup\left((U\times\R^m)\times\left(V\!\setminus\!(\supp v)\times\R^n\right) \right)\\
			&\qquad\cup \left((\supp u\times\R^m)\times (T^\ast V\!\setminus\{0\})\setminus\WF_{[\fM]}v\right)\\
			&\qquad\cup \left((T^\ast V\setminus\{0\})\setminus\WF_{[\fM]}u\times(\supp v\times\R^n)\right)
		\end{split}
	\end{equation*}
	then $(x_0,\xi_0,y_0,\eta_0)\notin\WF_{[\fM]}(u\otimes v)$.
	
	If $x_0\notin \supp u$ then there is a neighborhood $U^\prime\subseteq U$ of $x_0$
	such that $u\vert_{U^\prime}\equiv 0$. It follows that 
	\begin{equation*}
		(u\otimes v)\vert_{U^\prime\times V}\equiv 0
	\end{equation*}
	and thus $(x_0,\xi_0,y_0,\eta_0)\notin \WF_{[\fM]}(u\otimes v)$.
	Analogously we see that if $y_0\notin\supp v$ then $(x_0,\xi_0,y_0,\eta_0)\notin\WF_{[\fM]}(u\otimes v)$.
	
	Now, assume that $(x_0,\xi_0)\notin \WF_{(\fM)} u$. % and $\eta_0\neq 0$. 
	Then there are a neighborhood $U_1$ of $x_0$, 
	a conic neighborhood $\Gamma^\prime$
	of $\xi_0$ and a bounded sequence $(u_j)_j\subseteq\E^\prime(U)$
	with $u_j\vert_{U_1}=u\vert_{U_1}$ such that for all $\bM\in\fM$ and every $h>0$ there is a
	constant $C>0$ such that
	\begin{equation*}
		\xit^j\abs*{\hat{u}_j(\xi)}\leq Ch^jM_j
	\end{equation*}
	for all $\xi\in\Gamma^\prime$ and $j\in\Z_+$.
	
	Now choose a test function $\psi\in\CC^\infty_0(V)$ such that $\psi\vert_{V_1}\equiv 1$ for some
	neighborhood of $y_0$.
	Then $v_0=\psi v\in\E^\prime(V)$ and $w_j=u_j\otimes v_0\in\E^\prime(U\times V)$ is a bounded sequence with 
	\begin{equation*}
		w_j\vert_{U_1\times V_1}=(u\otimes v)\vert_{U_1\times V_1}.
	\end{equation*}
	Furthermore, we know that there are some $d\in\Z_+$ and $C_0>0$ such that
	\begin{equation*}
		\abs*{\hat{v}_0(\eta)}\leq C_0(1+\etat)^d\qquad \fa\eta\in\R^n.
	\end{equation*}
	For the tensor product we have that $\hat{w}_j(\xi,\eta)=\hat{u}_j(\xi)\hat{v}_0(\eta)$.
	We continue by choosing a conic neighborhood $\Gamma\subseteq\R^{m+n}$ of $(\xi_0,\eta_0)$
	such that $\pi_1(\Gamma)\subseteq\Gamma^\prime\subseteq\R^m$. Then it follows that
	\begin{equation*}
		\begin{split}
			\abs*{(\xi,\eta)}^j \abs*{\hat{w}_j(\xi,\eta)}&=\abs*{(\xi,\eta)}^j \abs*{\hat{u}_j(\xi)}
			\abs*{\hat{v}_j(\eta)}\\
			&\leq C_0Ch^{j+d} M_{j+d}\sup_{\substack{\xi^\prime\in\pi_1(\Gamma)\\ \eta^\prime\in\pi_2(\Gamma)}}
			\frac{\abs*{(\xi^\prime,\eta^\prime)}^j(1+\abs{\eta^\prime})^d}{\abs{\xi^\prime}^{j+d}}
		\end{split}
	\end{equation*}
	for $(\xi,\eta)\in\Gamma$. Now we set $\Gamma_1=\Set{\xi/\xit\given \xi\in\pi_1(\Gamma)}
	\subseteq S^{m-1}$
	and $\Gamma_2=\Set{\eta/\etat\given \eta\in\pi_2(\Gamma)}\subseteq S^{n-1}$.
	Then
	\begin{equation*}
		\begin{split}
			\sup_{\substack{\xi^\prime\in\pi_1(\Gamma)\\ \eta^\prime\in\pi_2(\Gamma)}}
			\frac{\abs*{(\xi^\prime,\eta^\prime)}^j(1+\abs{\eta^\prime})^d}{\abs{\xi^\prime}^{j+d}}
			&\leq \sup_{\substack{\xi^\prime\in\Gamma_1\\ \eta^\prime\in\Gamma_2\\\rho>1}}
			\frac{\abs*{(\rho\xi^\prime,\rho\eta^\prime}^j(1+\abs*{\rho\eta^\prime})^d}{\abs{\rho\xi^\prime}^{j+d}}\\
			&\leq \rho^{j+d-j-d}
			\sup_{\substack{\xi^\prime\in\Gamma_1\\ \eta^\prime\in\Gamma_2}}
			\abs*{(\xi^\prime,\eta^\prime)}^j\leq \tau^j
		\end{split}
	\end{equation*}
	for some $\tau>0$.
	Since $d$ only depends on $v$ we can apply $\eqref{B-DerivClosed}$ and conclude that
	for each $\bN\in\fM$ there are some $\bM\in\fM$ and
	a constant $\gamma>0$ such that $M_{j+d}\leq \gamma^{j+d+1}N_j$ for all $j\in\Z_+$. 
	Therefore we have shown that there is a conic neighborhood $\Gamma^\prime$ of $(\xi_0,\eta_0)$ 
	such that for all $\bN\in\fM$ and all $h>0$ we have that
	\begin{equation*}
		\sup_{\substack{(\xi,\eta)\in\Gamma^\prime\\ j\in\Z_+}}
		\frac{\abs*{(\xi,\eta)}^j\abs*{\hat{w}_j(\xi,\eta)}}{h^jN_j}<\infty.
	\end{equation*}
	Hence $(x_0,\eta_0,y_0,\eta_0)\notin\WF_{(\fM)} (u\otimes v)$.
	The Roumieu case follows similarly.
	
	If $(y_0,\eta_0)\notin\WF_{[\fM]}u$ and $\xi_0\neq 0$ we can argue analogously.
\end{proof}
For later reference we recall from \cite{MR4002151} some facts.
To begin with let
\begin{align*}
	I(\xi)&=\int_{\lvert\omega\rvert=1}\! e^{\omega\xi}\,d\omega,\\
	\tilde{K}(z)&=(2\pi)^{-n}\int\! e^{iz\xi}/I(\xi)\,d\xi.
\end{align*}
\begin{Thm}[{\cite[Theorem 5.7]{MR4002151}}]\label{BasicDecomposition} \label{WFThmChar1}
	If $u\in\mathcal{S}^\prime(\R^n)$ and $U=\tilde{K}\ast u$, then $U$ is analytic in 
	$X=\Set{z\in\C\given \abs{\imag z}<1}$ and there exist $C,a,b$ such that
	\begin{equation}\label{convest}
		| U(z)|\leq C\bigl(1+| z|\bigr)^{a}\bigl(1-|\imag z|)^{-b},\quad z\in X.
	\end{equation} 
	The boundary values $U(\,\cdot\,+i\omega)$ are continuous functions of $\omega\in S^{n-1}$ with values in $\mathcal{S}^\prime(\R^n)$, and
	\begin{equation}\label{DefConv}
		\langle u,\varphi\rangle=\int_{S^{n-1}}\!\langle U(\cdot+i\omega),\varphi\rangle\,d\omega,\quad \varphi\in \mathcal{S}.
	\end{equation}
	On the other hand, if $U$ is given satisfying \eqref{convest}, then the formula \eqref{DefConv} defines a distribution $u\in\mathcal{S}^\prime$ with $U=K\ast u$.
	
	For all
	$[$semiregular$]$ weight matrices $\fM$ we have
	\begin{equation*}
		\bigl(\R^n\times S^{n-1}\bigr)\cap \WF_{[\fM]} u= \bigl\{(x,\omega) : |\omega|=1,\, 
		U\text{ is not in }\E^{[\fM]} \text{ at } x-i\omega\bigr\}.
	\end{equation*}
\end{Thm}
Here $U$ being in $\E^{[\fM]}$ at $x-i\omega$ means that there is a neighborhood $V$ of $x-i\omega$ in $\C^n$ such that the derivatives $\partial_z^\alpha U$ satisfy
the defining estimates of $\E^{[\fM]}$ in $V$.
Theorem \ref{BasicDecomposition} follows from a straightforward modification of the proof of \cite[Theorem 8.4.11]{MR1996773}
using the [semiregularity] of $\fM$.
The same applies to the following corollary (cf.~\cite[Corollary 8.4.13]{MR1996773}).

\begin{Cor}[{\cite[Corollary 5.8]{MR4002151}}]\label{DecompCor}
	Let $\Gamma^1,\dotsc\Gamma^N\subseteq\R^n\!\setminus\!\{0\}$ be closed cones such that $\bigcup_j\Gamma^j=\R^n\!\setminus\!\{0\}$.
	Any $u\in\mathcal{S}^\prime(\R^n)$ can be written $u=\sum u_j$, where $u_j\in\mathcal{S}^\prime$ and
	\begin{equation}\label{WFdecomp}
		\WF_{[\fM]}u_j\subseteq\WF_{[\fM]}u\cap\bigl(\R^n\times\Gamma^j\bigr).
	\end{equation}
	If $u=\sum u^\prime_j$ is another such decomposition, then $u^\prime_j=u_j+\sum_k u_{jk}$, where $u_{jk}\in\mathcal{S}^\prime$, $u_{jk}=-u_{kj}$ and
	\begin{equation*}
		\WF_{[\fM]}u_{jk}\subseteq \WF_{[\fM]}u\cap\bigl(\R^n\times\bigl(\Gamma^j\cap\Gamma^k\bigr)\bigr).
	\end{equation*}
\end{Cor}
The next theorem generalizes \cite[Theorem 8.4.15]{MR1996773}; it suffices to follow the 
arguments in \cite{MR1996773} while taking the above statements into account;
 recall that $\Gamma^\circ := \Set{\xi \in \R^n \given \langle y, \xi \rangle \ge 0 \text{ for all } y \in \Gamma}$ 
denotes the dual cone of $\Gamma$.
Moreover a holomorphic function $F\in \mathcal{O}(V+i\Gamma_r)$ is said to be of slow growth if 
\begin{equation*}
	\ex C, \gamma>0:\quad \abs{F(x+iy)}\leq C\abs{y}^{-\gamma},\quad y\in\Gamma_r=\{\tilde{y}\in\Gamma\given\abs{\tilde{y}}<r\}.
\end{equation*}
\begin{Cor}\label{SemiBVChar} 
	Let $\fM$ be a [semiregular] weight matrix.
	Let $u\in\D^\prime(\Omega)$ and $(x_0,\xi_0)\in T^\ast\Omega \setminus \{0\}$.
	Then $(x_0,\xi_0)\notin\WF_{[\fM]} u$ if and only if there exists a neighborhood $V$ of $x_0$, $f\in\E^{[\fM]}(V)$,
	open cones $\Gamma^1,\dotsc,\Gamma^q$ with the property $\xi_0\Gamma^j<0$ for all $j$, and holomorphic functions
	$F_j\in\mathcal{O}(\Omega+i\Gamma^j_{\delta})$ of slow growth such that
	\begin{equation*}
		u\vert_V=f+\sum_{j=1}^qb_{\Gamma^j}F_j.
	\end{equation*}
\end{Cor}

In \cite{MR1996773}the pullback of distributions with real analytic mappings was defined.
Since $\E^{[\fM]}$ is stable by pullback with real analytic mappings, see \cite[Theorem 2.9]{MR4002151}, 
we can follow the proof of \cite[Theorem 8.5.1]{MR1996773} to obtain the following statement:

\begin{Thm}[{\cite[Theorem 5.11]{MR4002151}}]\label{Invar1}
	Let $\fM$ be a $[$semiregular$]$ weight matrix.
	Let $F:\Omega_1\rightarrow\Omega_2$ be a real analytic mapping, where $\Omega_i \in \R^{n_i}$ are open.
	If $u\in\D^\prime(\Omega_2)$ and $N_F\cap\WF_{[\fM]}u=\emptyset$, then	
	\begin{equation*}
		\WF_{[\fM]}\left(F^\ast u\right)\subseteq F^\ast\left( \WF_{[\fM]}u\right). 
	\end{equation*}
	Here $N_F = \Set{ (f(x),\eta) \in \Omega_2 \times \R^{n_2}\given F'(x)^T \eta = 0 }$ is the set of normals of $F$.
\end{Thm}
\begin{Thm}\label{WFproduct}
	Let $\fM$ be a [semiregular] weight matrix.
	If $u,v\in\D^\prime(\Omega)$ are such that $(x,\xi)\in\WF_{[\fM]}u$ implies that $(x,-\xi)\notin\WF_{[\fM]}v$ then the product $uv$ is well defined as a distribution in 
	$D^\prime(\Omega)$ and
	\begin{equation*}
		\WF_{[\fM]} (uv)\subseteq\Set*{(x,\xi+\eta)\in\Omega\times\R^n\!\setminus\!\{0\}
			\given (x,\xi)\in\WF_{[\fM]}u\text{ or }\xi=0,\; (x,\eta)\in\WF_{[\fM]}v\text{ or }\eta=0}.
	\end{equation*}
\end{Thm}
\begin{proof}
	 \cite[Theorem 8.2.10]{MR1996773} gives that $uv$ is well-defined.
	There the product $uv$ is defined as
	\begin{equation*}
		\delta^\ast(u\otimes v) 
	\end{equation*}
	where $\delta(x)=(x,x)$ is the diagonal map from $\Omega$ to $\Omega\times\Omega$.
	The statement follows from the proof of \cite[Theorem 8.2.10]{MR1996773}
	when we apply Theorem \ref{WFTensor} and Theorem \ref{Invar1}.
\end{proof}
\begin{Prop}\label{BasicWFIntegral}
Let $u\in\E^\prime(\R^{n+m})$ and denote the variables of $\R^{n+m}$ by $(x,y)$. We define 
\begin{equation*}
	u_1=\int_{\R^m}\!u(\:.\:,y)\,dy\in\E^\prime(\R^n).
\end{equation*}
	(More precisely, $\langle u_1,\phi\rangle=\langle u,\phi\otimes 1\rangle$).
	Then 
	\begin{equation*}
		\WF_{[\fM]}u_1=\Set*{(x,\xi)\in\R^n\times\R^n\!\setminus\!\{0\} 
		\given (x,y,\xi,0)\in\WF_{[\fM]}u\text{ for some }y\in\R^m}
	\end{equation*}
	for all [semiregular] weight matrices $\fM$.
\end{Prop}
\begin{proof}
	Using the notation of \cite[Theorem 5.7]{MR4002151} we have for $\phi\in\CC^\infty_0(\R^n)$ and $\psi\in\CC^\infty_0(\R^m)$ that
	\begin{equation*}
		\langle u,\phi\otimes\psi\rangle=\int_{\lvert\omega\rvert=1}\!\langle U(\,.\,+i\omega),\phi\otimes\psi\rangle\,d\omega.
	\end{equation*}
	Take $\psi(y)=\chi(\delta y)$ where $\chi$ is a test function that equals $1$ on the unit ball, and let $\delta\rightarrow 1$.
	Since $U$ is exponentially decreasing at infinity it follows that
	\begin{equation*}
		\langle u_1,\phi\rangle =\int_{\lvert\omega\rvert=1}\!\langle U(\,.\,+i\omega),\phi\otimes 1\rangle\,d\omega
	=\int_{\lvert\tau\rvert=1}\! U_1(\, .\,+i\sigma),\phi\rangle\,d\sigma
\end{equation*}
	where
	\begin{equation*}
		U_1(\zeta)=\int\! U(\zeta,y)\,dy=\int\! U(\zeta,y+\tau)\,dy,\qquad\quad \lvert\imag\zeta\rvert^2+\lvert\tau\rvert^2<1,
	\end{equation*}
	is an analytic function when $\lvert\imag\zeta\rvert<1$ which is bounded by $C(1-\lvert\imag\zeta\rvert)^{-N}$.
	If $\lvert\sigma\rvert=1$ and $(x,y,\sigma,0)\notin\WF_{[\fM]} u$ for all $y\in\R^m$ then $U_1\in\E_{[\fM]}$ at $x-i\sigma$
	by Theorem \ref{BasicDecomposition}.
	Hence the straightforward generalization of \cite[Lemma 8.4.12]{MR1996773} implies $(x,\sigma)\notin\WF_{[\fM]}u_1$.
 \end{proof}

For the next statement we need more definitions: Let $U\subseteq\R^m$ and $V\in\R^n$ be open
sets and $K\in\D^\prime(U\times V)$. We set
\begin{align*}
	\WF_{[\fM]}(K)_U&:=\Set*{(x,\xi)\in U\times\R^m\!\setminus\!\{0\}
		\given (x,y,\xi,0)\in\WF_{[\fM]} K \text{ for some }y\in V}\\
	\WF_{[\fM]}^\prime(K)_V&:=\Set*{(y,\eta)\in V\times\R^n\!\setminus\!\{0\}
		\given (x,y,0,-\eta)\in \WF_{[\fM]} K\text{ for some }x\in U}\\
	\WF_{[\fM]}^\prime(K)&:=\Set*{(x,y,\xi,\eta)\in T^\ast(U\times V)\setminus\{0\}\given
		(x,y,\xi,-\eta)\in\WF_{[\fM]}(K)}
\end{align*}
It is clear that $K$ defines a mapping $\mathcal{K}:\CC^\infty_0(V)\rightarrow \D^\prime(U)$.
Unter certain assumptions on the wavefront sets of $K$ and $u\in\D^\prime(V)$ it is possible to extend
$\mathcal{K}$ to make sense of $\mathcal{K}u$, cf.~\cite[Theorem 8.2.13 \& Theorem 8.5.5]{MR1996773}.
The following Theorem is the main statement of this section and is a direct generalization of \cite[Theorem 8.5.5]{MR1996773}.
\begin{Thm}\label{WFkernel}
	If $u\in\E^\prime(V)$ and $\WF_{[\fM]}u\cap\WF_{[\fM]}^\prime(K)_V=\emptyset$
	then $\mathcal{K}u\in\D^\prime(U)$ is well-defined
	\begin{equation*}
		\begin{split}
			\WF_{[\fM]}(\mathcal{K}u)&\subseteq\WF_{[\fM]}(K)_U\cup\left(
			\WF_{[\fM]}^\prime(K)\circ\WF_{[\fM]}u\right)\\
			&=\WF_{[\fM]}(K)_U\cup\Bigl\{(x,\xi)\in U\times\R^n\setminus\{0\}\;:\\
			&\qquad\qquad\qquad	(x,y,\xi,-\eta)\in\WF_{[\fM]}(K)\text{ for some }(x,\eta)\in\WF_{[\fM]}u\Bigr\}.
		\end{split}
	\end{equation*}
\end{Thm}
\begin{proof}
	By Theorem \ref{WFTensor} we have that
	\begin{equation*}
		\WF_{[\fM]}(1\otimes u)\subseteq\Set*{(x,y,0,\eta)\given (y,\eta)\in\WF_{[\fM]}u}.
	\end{equation*}
	Thence $K_u=K(1\otimes u)$ is well-defined and Theorem \ref{WFproduct} gives 
	\begin{equation*}
		\begin{split}
		\WF_{[\fM]}(K_u)&\subseteq\Set*{(x,y,\xi,\eta+\eta^\prime)\given
			(y,\eta)\in\WF_{[\fM]}u,\;(x,y,\xi,\eta^\prime)\in\WF_{[\fM]}K}\\
			&\qquad
		\cup\WF_{[\fM]}K\cup\WF_{[\fM]}(1\otimes u).
	\end{split}
	\end{equation*}
	Now $\langle \mathcal{K}u,\psi\rangle=\langle K_u,\psi\otimes 1\rangle$ by the definition of 
	$\mathcal{K}$ in \cite[Theorem 8.2.13]{MR1996773}.
	Therefore the assertion follows from Proposition \ref{BasicWFIntegral}.
\end{proof}

\begin{Rem}
	If $\fM=\{\bM\}$ consists of a single weight sequence then
	in the Roumieu case we just recover the results of \cite[Section 8.5]{MR1996773}.
	In the case of weight functions we need a little more work, 
	cf.~\cite{MR3285413}:
	
	If $\omega$ is a weight function then we associate to $\omega$ a weight matrix 
	$\mathfrak{W}=\Set{\mathbf{W}^\lambda=(W_j^\lambda)_j\given \lambda>0}$ by setting
	\begin{equation*}
		W_j^\lambda =\exp\left(\frac{1}{\lambda}\varphi_\omega^\ast(\lambda j)\right).
	\end{equation*}
	Then $\mathfrak{W}$ satisfies the following properties
	\begin{gather*}
		\fa\lambda>0\,\fa j_1,j_2\in\Z_+\quad W^{\lambda}_{j_1+j_2}\leq W_{j_1}^{2\lambda}
		W_{j_2}^{2\lambda}\\
		\fa h\geq 1 \,\ex A\geq 1\, \fa\lambda>0\, \ex D\geq 1 \,\fa j\in\Z_+:\quad
		h^jW_j^\lambda\leq DW_j^{A\lambda}.
	\end{gather*}
	Thus $\DC{\omega}{\Omega}$ and $\DC{\mathfrak{W}}{\Omega}$ are isomorphic as topological vector spaces.
	Furthermore $\mathfrak{W}$ is $R$- and $B$-semiregular
	 if $\omega(t)=o(t)$, $t\rightarrow\infty$.
	Finally it is easy to see that $\WF_{[\mathfrak{W}]}u=\WF_{[\omega]}u$ where
	$\WF_{[\omega]}u$ is the wavefront set associated to $\DC{\omega}{\Omega}$ defined in
	\cite{Albanese2010}.
\end{Rem}

We close this section by presenting a few applications of the main Theorem \ref{WFkernel}.
First, we observe that, by using Theorem \ref{WFkernel} and Theorem \ref{Invar1}, we can easily extend the arguments 
from  \cite[p.~270]{MR1996773} and obtain the following statement regarding the convolution of distributions.

\begin{Cor}\label{ConvWF}
	Let $k\in\D^\prime(\R^n)$ and $u\in\E^\prime(\R^n)$.
	Then
	\begin{equation*}
		\WF_{[\fM]}(k\ast u)\subseteq\Set*{(x+y,\xi)\given (x,\xi)\in\WF_{[\fM]} k\text{ and }
			(y,\xi)\in\WF_{[\fM]}u}
	\end{equation*}
\end{Cor}
We apply Corollary \ref{ConvWF} to the asymptotically analytic cutoff multipliers from 
\cite[Section 17.3]{MR4436039}.
More precisely, given two open cones $\Gamma,\Gamma^\prime\subseteq\R^n\setminus\{0\}$ with
$\Gamma^\prime\cap S^{n-1}\Subset\Gamma\cap S^{n-1}\neq S^{n-1}$ (We assume here that $n\geq 2$) and
$R>0$ there is a smooth function  $g_R\in\CC^\infty(\C\!\setminus\! i\R^n)$ with the following properties (cf.~\cite[Lemma 17.3.2]{MR4436039}):
\begin{itemize}
	\item $\fa\xi\in\R^n\setminus\{0\}:\quad 0\leq g_R(\xi)\leq 1$
	\item  $\fa\xi\in\R^n\setminus\{0\},\fa \alpha\in\Z_+^n,\;\xit\geq 2R\alp:\quad
	\abs{D^\alpha g_R(\xi)}\leq 2(C_1/R)^{\alp}$
	\item $\fa\zeta=\xi+i\eta\in\C^n\!\setminus\!i\R^n$:
	\begin{itemize}
		\item If $\xi\in\Gamma^\prime$ then $g_R(\xi)=1$.
		\item If $\xi\notin\Gamma$ then $g_R(\xi)=0$.
	\end{itemize}
	\item $\fa\zeta=\xi+i\eta\in\C^n\!\setminus\!i\R^n$:
	\begin{gather*}
		\abs*{g_R(\zeta)}\leq \exp\left(C_1\frac{\abs{\eta}}{R}\right)\\
		\abs*{\partial_{\zeta}g_R(\zeta)}\leq C_2\exp\left((C_3\abs{\eta}-\abs{\xi})/R\right)
	\end{gather*}
\end{itemize}
The constants $C_1$, $C_2$ and $C_3$ do not depend on $R$.

If we denote by $G_R=\mathcal{F}^{-1}(g_R)$ the inverse Fourier transform of $g_R$ then
we obtain from \cite[Lemma 17.3.6]{MR4436039} that $\WF_a G_R\subseteq \R^n\times\overline{\Gamma}$ since $G_R\ast\delta=G_R$.
Moreover, $G_R$ is analytic in the region $\abs{x}> C/R$ where $C$ is a constant independent of $R$ cf.~\cite[Proposition 17.3.3]{MR4436039}. We obtain that 
\begin{equation*}
	\WF_{[\fM]} G_R\subseteq \Set*{x\in\R^n\given \abs{x}\leq C/R}%\given \overline{\Gamma}
\end{equation*}
for any [semiregular] weight matrix.
Following the presentation in \cite[Section 17.3]{MR4436039} we write
\begin{equation*}
	g_R(D)u=G_R\ast u,\qquad u\in\E^\prime(\R^n).
\end{equation*}
\begin{Prop}\label{Multiplier}
	Let $u\in\E^\prime(\R^n)$, $(x_0,\xi_0)\in \R^n \times\R^n\setminus\{0\}$ and $\fM$
	be a [semiregular] weight sequence.
	Then $(x_0,\xi_0)\notin\WF_{[\fM]}u$ if and only if there are an open neighborhood of
	$x_0$, two open cones $\xi_0\in\Gamma^\prime\Subset\Gamma\subsetneq\R^n\setminus\{0\}$ 
	and some $R>0$ such that $g_R(D)u\vert_U\in\DC{\fM}{U}$.
\end{Prop}
\begin{proof}
	If $(x_0,\xi_0)\notin\WF_{[\fM]} u$ then there is an open neighborhood $V$ of $x_0$
	and some open convex cone $\xi_0\in F\subsetneq\R^n\setminus\{0\}$ such that
	$\WF_{[\fM]}u\cap V\times F\neq\emptyset$.
	We choose two open convex cones $\Gamma,\Gamma^\prime$ such that
	$\xi_0 \in\Gamma^\prime\Subset\Gamma\Subset F$. 
	For $R$ to be determined later we take the multiplier $g_R$ from above which is 
	associated to the cones $\Gamma$ and $\Gamma^\prime$.
	We have shown that $\WF_{[\fM]} G_R\subseteq\Set{(x,\xi)\given \abs{x}\leq C/R,\; \xi\in\overline{\Gamma}}$.
	If we choose $R>0$ such that $C/R<\tfrac{1}{2}\dist (x_0,\R^n\setminus V)$ then
	it follows that there is a neighborhood $W$ of $x_0$ such that 
	\begin{equation*}
		\WF_{[\fM]}\bigl( G_R\ast u\bigr)\cap W\times\R^n\setminus\{0\}=\emptyset.
	\end{equation*}

The other direction can be proven in the same way as the sufficiency part of the proof of
\cite[Proposition 17.3.7]{MR4436039}: If $g_R$ is an multiplier of the above form for the cones $\Gamma^\prime\Subset\Gamma$, then the same is true for $1-g_R$
with respect to the cones $\emptyset\neq\R^n\setminus\Gamma\Subset\R^n\setminus\Gamma^\prime\subsetneq\R^n\setminus\{0\}$.
Thence
\begin{align*}
\WF_{[\fM]} (1-g_R(D))u\cap \R^n\times\Gamma^\prime&=\emptyset\\
\shortintertext{and thus}
\WF_{[\fM]} u\cap \R^n\times\Gamma^\prime&=\emptyset.
\end{align*} 
\end{proof}

Finally we are now able to prove Theorem \ref{Cordaro1} which will follow from the combination of
\cite[Theorem 5]{MR2126468} and the following statement:
\begin{Thm}
	Let $k\geq 2$ be an even integer and $\fM$ be a [semiregular] weight matrix.
	If $p_0=(x_0,0;\xi_0,0)\in T^\ast\R^2\!\setminus\!\{0\}$
	then the following statements are equivalent:
	\begin{enumerate}
		\item $Q_\lambda$ is $[\fM]$-hypoelliptic at $p_0$, i.e.~for any distribution $u$
		defined near $(x_0,0)$ we have that $p\notin\WF_{[\fM]} Q_\lambda u$ implies
		$p\notin\WF_{[\fM]} u$
		\item $\lambda\notin\Lambda_k$
	\end{enumerate}
\end{Thm}
\begin{proof}
	Given $\lambda\in\Lambda_k$ and $p_0=(x_0,0;\xi_0,0)$ then according to \cite[Lemma 1]{MR2126468} there exists a distribution $u\in\D^\prime(\R^n)$ such that $Q_\lambda u=0$
	and 
	\begin{align*}
		\WF u=\WF_a u&=\Set*{(x_0,0;\rho\xi_0,0)\given \rho>0}.\\
		\shortintertext{Thus Proposition \ref{WFproperties}(4) gives that}
		\WF_{[\fM]}u&=\Set*{(x_0,0;\rho\xi_0,0)\given \rho>0}
	\end{align*}
	for any [semiregular] weight matrix $\fM$. Therefore the implication $(1)\Rightarrow(2)$ is proven.
	If $(2)$ holds then in \cite[Section 5]{MR2126468} a left parametrix for $Q_\lambda$ is constructed: 
	\begin{equation*}
		\cK Q_\lambda=\Id +\cR
	\end{equation*}
	in $\mathcal{S}(\R^2)$ and the kernels $K$ and $R$ of $\cK$ and $\cR$, respectively, satisfy
	\begin{align*}
		\WF K=\WF_aK&=\Set*{(x,t,x^\prime,t^\prime;\xi,\tau,\xi^\prime,\tau^\prime)\given
			x=x^\prime,\; t=t^\prime,\; \xi=-\xi^\prime,\; \tau=-\tau^\prime},\\
		\WF R=\WF_a R&=\Set*{(x,t,x^\prime,t^\prime;0,\tau,0,\tau^\prime)\given x=x^\prime,\; t=t^\prime,\; \tau=-\tau^\prime},
	\end{align*}
	see \cite[equations (33) \& (35)]{MR2126468}. Now (1) follows from Theorem \ref{WFproperties}(3)
	and Theorem \ref{WFkernel}.
\end{proof}
\section{Analytic pseudodifferential operators and ultradifferentiable microdistributions}\label{Sec:MicroDistr}
An obvious application of Theorem \ref{WFkernel} would be,
for example, that $s$-Gevrey pseudodifferential operators 
(for a definition and details see e.g. \cite{MR1249275})  is
$[\fM]$-pseudolocal with respect to 
any ultradifferentiable class $\mathcal{E}^{[\fM]}\supseteq\G^s$, which is closed
under differentiation.
The next step would be to consider ultradifferentiable
pseudodifferential operators but this would involve 
 more conditions on the data of the ultradifferentiable
classes than we considered in the section before. Thus
we relegate this topic to a forthcoming paper and focus in this paper, as noted in the introduction, on analytic operators, i.e.~generally operators with kernels which analytic wavefront sets are ``well-behaved''.
Such operators include in particular analytic differential operators. For completeness, we have included the following 
theorem  which gives a very large family of operators
$A$ such that $\WF_{[\fM]}Au\subseteq\WF_{[\fM]}u$
for every [semiregular] weight matrix $\fM$
by Theorem \ref{WFkernel} and Proposition \ref{WFproperties}(4).
\begin{Thm}\label{Osc}
	Let $\Omega\subseteq\R^n$ be an open set and
	 $\Phi\in\CC^\infty(\Omega\times\R^N\!\setminus\!\{0\})$ with $\imag \Phi\geq 0$.
	 Suppose that $\Phi$ depends holomorphically on the first variable, is homogeneous of degree $1$ in the second variable and $\nabla_x\Phi(x,\theta)\neq 0$ on
	 $\Delta_\Phi=\Set{(x,\theta)\given d_\theta\Phi(x,\theta)=0}$.
	 Furthermore let $a\in\CC^\infty(\Omega\times\R^N)$ be
	 such that $a$ extends to a holomorphic function in the first
	 variable on some complex neighborhood $\widetilde{\Omega}$ of $\Omega$
	 and $a$ satisfies the following estimates:
	 There is some $d\in\R$ such that
	 for each compact set $K\subseteq\widetilde{\Omega}$, 
	 and all
	 $\gamma\in\Z_+^N$ there is a constant $C_\gamma>0$ so that
	 \begin{equation}\label{GeneralEstimates}
	 	\abs*{\partial_\theta^\gamma a(z,\theta)}\leq C_\gamma
	 	\left(1+\abs{\theta}\right)^{d-\abs{\gamma}} \qquad\fa z\in K,\; \fa\theta\in\R^N.
	 \end{equation}
	 
	 Then the oscillatory integral 
	 \begin{gather*}
	 	A(x)=\int e^{i\Phi(x,\theta)}a(x,\theta)\,d\theta\in \D^\prime(\Omega)\\
	 	 \shortintertext{exists and}
	 	 \WF_a A\subseteq\Set*{(x,\real d_x\Phi(x,\theta))\in T^\ast\Omega\setminus\{0\}\given d_\theta\Phi(x,\theta)=0}.
	 \end{gather*}
	 \end{Thm}
	\begin{proof}
	First note that if $\psi\in\CC^\infty_0(\R^N)$
	with $\psi(\theta)=1$ for $\abs{\theta}< 1$
	then 
	\begin{equation*}
		h_0(x)=	\int \psi(\theta)e^{i\Phi(x,\theta)}a(x,\theta)\,d\theta
		\end{equation*}
		extends to a holomorphic function in a neighborhood of $\Omega$.
		Hence we can w.l.o.g.~assumme that $a(x,\theta)$
		vanishes for $\abs{\theta}<1$.
		
		Now assume that $x_0\in\Omega$ such that $d_\theta\Phi(x_0,\theta)\neq 0$ for all $\theta\in\R^N\!\setminus\!\{0\}$.
		That means that either $d_\theta\real\Phi(x_0,\theta)\neq 0$ or $d_\theta\imag\Phi(x_0,\theta)\neq 0$ for all $\theta$.
		In the latter case we have necessarily that $\imag \Phi(x_0,\theta)>0$ for all $\theta\in\R^N\!\setminus\!\{0\}$.
		Let $V\subseteq\widetilde{\Omega}$ be a 
		complex neighborhood of $x_0$ such that
		$\imag\Phi(z,\vartheta)>0$ for all $z\in \overline{V}$
		and $\vartheta\in S^{N-1}$. Then there exists some $c>0$
		such that
		\begin{align*}
			\abs*{\imag \Phi(z,\vartheta)}&>c & \fa \vartheta&\in S^{N-1},\;\fa z\in \overline{V}\\
			\shortintertext{and therefore}
			\abs*{\imag \Phi(z,\theta)}&>c\abs{\theta}
			& \fa \theta&\in\R^N\!\setminus\!\{0\},\;\fa z\in \overline{V} 
		\end{align*}
		by homogeneity.
		Thence
		\begin{equation*}
			A(z)=\int e^{i\Phi(z,\theta)}a(z,\theta)\,d\theta
		\end{equation*}
		is a holomorphic function in $V$.
		
		In the first case %we observe that there is some $k\in\{1,\dotsc,N\}$ such that
		%$\partial_{\theta_k}\real\Phi(x_0,\theta)$ for $\theta\in\R^N\!\setminus\!\{0\}$.
		we choose again a complex neighborhood $V\subseteq\widetilde{\Omega}$ of $x_0$
		such that
		\begin{equation}\label{LowerBound1}
			\abs*{\partial_{\theta}\Phi(z,\theta)+2i\eps\theta}> c \qquad \fa z\in\overline{V},\;\fa\theta\in S^{N-1},\;\fa \eps\in\R,
		\end{equation}
		for some constant $c>0$.
		For $\eps>0$ we set
		\begin{equation*}
		A_\eps(x)=\int\!e^{i\Phi(x,\theta)-\eps\abs{\theta}^2}
			a(x,\theta)\,d\theta.
		\end{equation*}
		By \cite[Theorem 8.1.13]{MR4436039} $A_\eps\rightarrow A$
		in $\D^\prime(V)$.
		Obviously $A_\eps$ extends to a holomorphic function
		in $\widetilde{\Omega}$. 
		We follow primarily the arguments presented in \cite[Section 18.1]{MR4436039}:
		By integration by parts we have that
		\begin{equation*}
			\begin{split}
				A_\eps(z)&=
				\sum_{j=1}^N\int D_{\theta_j}\left(e^{i\Phi(z,\theta)-\eps\abs{\theta}^2}\right)\frac{a(z,\theta)}{\sum_{k=1}^N\partial_{\theta_k}\Phi(z,\theta)+2i\eps\theta_k}\,d\theta\\
				&=-\sum_{j=1}\int e^{i\Phi(z,\theta)-\eps\abs{\theta}^2} 
				D_{\theta_j}\left(\frac{a(z,\theta)}{
					\sum_{k=1}^N\partial_{\theta_k}
				\Phi(z,\theta)+2i\eps\theta_k}\right)\,d\theta.
			\end{split}
		\end{equation*}
		Set $\Psi(z,\theta)=(\sum_{k=1}^N\partial_{\theta_k}
		\Phi(z,\theta)+2i\theta_k)^{-1}$.
		Applying \eqref{LowerBound1} we conclude that
		\begin{equation*}
			D^\alpha\Psi(z,\theta)=O(\abs{\theta}^{-\abs{\alpha}}),\qquad \abs{\theta}\longrightarrow\infty
		\end{equation*}
		Therefore,
		using also \eqref{GeneralEstimates}, we obtain that
		\begin{equation*}
			D_{\theta_k}\left(\frac{a(z,\theta)}{\partial_{\theta_k}
				\Phi(z,\theta)+2i\eps\theta_k}\right)
		\end{equation*}
		satisfies also estimates of the form \eqref{GeneralEstimates}
		with $d$ replaced by $d-1$ for all
		$z\in\overline{V}$ and every $\eps\in\R$.
		Iterating this procedure we obtain that
		\begin{equation*}
		A_\eps(z)=\int\!e^{i\Phi(z,\theta)-\eps\abs{\theta}^2}
		a_\eps(z,\theta)\,d\eps,\qquad \eps>0,
		\end{equation*}
		where $a_\eps$, $\eps\in\R$, are functions
		satisfying estimates of the form \eqref{GeneralEstimates}
		with $d$ replaced by $-N-1$ and 
		\begin{equation*}
			\sup_{z\in V}\abs{g_\eps(z,\theta)-g_0(z,\theta)}\leq
			C_\eps (1+\abs{\theta})^{-N-1},\qquad \theta\in\R^N\!\setminus\!\{0\}
		\end{equation*}
		where $C_\eps\rightarrow 0$ for $\eps\rightarrow 0$.
		Thus 
		\begin{equation*}
			A_\eps(z)\longrightarrow \int\! e^{i\Phi(z,\theta)}a_0(z,\theta)\,d\theta,\qquad
			\eps\rightarrow 0,
		\end{equation*}
		uniformly on $K$.
		Thence $\singsupp_a A\subseteq \pi_x(\Delta_\Phi)$.

It remains to analyze $A$ near a point $x_0\in\Omega$
for which there is some $\theta\in\R^N\!\setminus\!\{0\}$
with $d_\theta\Phi(x,\theta)=0$.
We choose a small neighborhood $V$ of $x_0$ and
 define 
\begin{equation*}
	\Delta_V=\Set*{\theta\in\R^N\!\setminus\!\{0\}\given
	\ex x\in V \text{ s.t. } d_\theta(x,\theta)=0}
\end{equation*}
which is obviously a closed cone in $\R^N\!\setminus\!\{0\}$.

We consider a covering of $\Delta_V$ by convex, open cones:
Let $F_1,\dotsc,F_\nu$ be open convex cones such that
\begin{gather*}
	\Delta_V\subseteq \bigcup_{j=1}^\nu F_j\\
	\shortintertext{and}
	d_x\Phi(x,\theta)\neq 0,\qquad \fa (x,\theta)\in V\times\bigcup_{j=1}^\nu F_j.
\end{gather*}
We choose smooth functions $\chi_j\in\CC^\infty(\R^N)$ such that
$0\leq \chi_j\leq 1$, $\chi_j(\theta)=0$ if $\abs{\theta}<1/2$
or if $\theta\notin\ \R^N\setminus F_j$ for all $1\leq j\leq\nu$
and
\begin{equation*}
	\sum_{j=1}^\nu\chi_j(\theta)=1,\qquad \fa \theta\in\Delta_V,\;\abs{\theta}\geq 1.
\end{equation*}
We set $\gamma=1-\sum\chi_j$. Thus formally
\begin{equation*}
	\int\!e^{i\Phi(x,\theta)}a(x,\theta)\,d\theta=
	\int\!e^{i\Phi(x,\theta)}\gamma(\theta)a(x,\theta)\,d\theta
	+\sum_{j=1}^\nu\int\!e^{i\Phi(x,\delta)}\chi_j(\theta)
	a(x,\theta)\,d\theta.
\end{equation*}
Following the arguments in the first part of the proof
 we can show that the first integral on the right-hand side defines a holomorphic function
in a complex neighborhood of $x_0$.

It remains to look at the oscillatory integral
 $B_j=\int e^{i\Phi}b_j(x,\theta)\,d\theta$ where we have set
 $b_j(x,\theta)=\chi_j(\theta)a(x,\theta)$ for $j=1,\dotsc,\nu$.
We begin by considering the Taylor expansion of $\Phi(x+iy,\theta)$ with respect to $y$:
\begin{equation*}
	\Phi(x+iy,\theta)=\Phi(x,\theta)+\nabla_y\Phi(x,\theta)y+ O(\abs{y}^2),\qquad x+iy\in\widetilde{\Omega},\theta\in\R^N\setminus\{0\}.
\end{equation*}
Since $\Phi(x+iy,\theta)$ is holomorphic in the first variable 
we have that $\nabla_y\Phi=i\nabla_x\Phi$. Thence
\begin{equation*}
	\imag \Phi(x+iy,\theta)=\imag \Phi(x,\theta)+\nabla_x\real\Phi(x,\theta)y+O(\abs{y}^2)
\end{equation*}
It follows that 
\begin{equation*}
	G_j(x+iy)=\int e^{i\Phi(x+iy,\theta)}b_j(x+iy,\theta)\,d\theta
\end{equation*}
is holomorphic in the set $\Omega+ i\Gamma^j_\delta$ for $\delta>0$
small enough. Here 
\begin{equation*}
	\Gamma^j=\Set*{y\in \R^n\given y\cdot v>0\;\fa v\in V_j}
\end{equation*}
with $V_j=\Set{\real\nabla_x\Phi(x,\theta)\given (x,\theta)\in V\times F_j}$. Obviously $V_j$ and $\Gamma^j$ are open cones in
$\R^n\!\setminus\!\{0\}$ and we write
also $\Gamma^j_\delta=\Set{y\in\Gamma\given \abs{y}<\delta}$.

We denote the boundary value (as a hyperfunction in $V$) of $G_j$ by $\bv G_j$.
Recall from \cite[Chapter 17]{MR4436039} that a representative $\mu$
of the restriction $\bv G_j$ to any open subset 
$V^\prime\subseteq V$ is given by the limit 
 in $\mathcal{O}^\prime(\C^n)$ of
the analytic functionals
\begin{equation*}
	\mathcal{O}(\C^n)\ni \varphi\xlongrightarrow{\quad\mu_t\quad}\int_{V^\prime+itv}G_j(w)\varphi(w)\,dw,\qquad (1>t\rightarrow +0).
\end{equation*}
Now observe that,
if we denote
\begin{equation*}
	G_j^\eps(z)=\int\!e^{i\Phi(z,\theta)-\eps\abs{\theta}^2}
	b_j(z,\theta)\,d\theta,
\end{equation*}
then $G^\eps_j$ converges to $G_j$ in $\mathcal{O}(V+i\Gamma^j_\delta)$.
Hence the analytic functionals 
\begin{equation*}
	\mu_t^\eps(\varphi)=\int_{V^\prime+itv}\! G^j(w)\varphi(w)\,dw,\qquad \varphi \in\mathcal{O}(\C^n),
\end{equation*}
converge to $\mu$ in $\mathcal{O}(\C^n)$ for $t,\eps\rightarrow\infty$.
On the other hand, the holomorphic functions $G_j^\eps(x+iy)$ converge
for $y\rightarrow 0$ to the smooth function
\begin{equation*}
	B_j^\eps(x)=\int\!e^{i\Phi(x,\theta)-\eps\abs{\theta}^2}b_j(x,\theta)\,d\theta
\end{equation*}
in $\CC(V)$ and $B_j^\eps$ converge to $B_j$ in $\D^\prime(V)$.
Now, choose for any $V^\prime\Subset V$ a function $\psi\in\CC^\infty_0(V)$ with $\psi\vert_{V^\prime}=1$.
It follows that $\psi B_j^\eps(\,.\,+iy)\rightarrow \psi B$
in $\E^\prime(\R^n)$ for $y,\eps\rightarrow 0$.
By the canonical embedding of $\E^\prime(\R^n)$ into
$\mathcal{O}^\prime(\C^n)$ we have that $\psi B_j^\eps(\,.\,+iy)$ converges to $\psi B$ in $\mathcal{O}^\prime(\C^n)$.
Therefore $\bv G_j$ coincide with $B_j$ as hyperfunctions in some
neighborhood of $x_0$. 
When we may shrink $V$ if necessary, we conclude that
$B_j\vert_V$ is the distributional boundary value of $G_j$
(cf.~\cite[Theorem 7.4.15]{MR4436039}) and
by \cite[Theorem 3.5.5]{MR4436039} we have that
\begin{equation*}
	\WF_a B_j\subseteq V\times (\Gamma^j)^\circ
\end{equation*}
where $(\Gamma^j)^\circ$ is the polar cone of $\Gamma^j$, i.e.
\begin{equation*}
	(\Gamma^j)^\circ=\Set*{v\in\R^n\given yv\geq 0\;\fa y\in\Gamma^j}=\overline{V_j}.
\end{equation*}

In order to finish the proof note that, if $\xi_0\in\R^n\!\setminus\{0\}$ is such that $(x_0,\xi_0)\notin\Delta_\Phi$ then
we can choose the cones $F_j$, $j=1,\dotsc,\nu$, in such a way 
that $\xi_0\notin \overline{V_j}$ for all $j=1,\dotsc,\nu$.
Thence $(x_0,\xi_0)\notin\WF_a A$.
\end{proof}

The class of distributions appearing in Theorem \ref{Osc}
includes the distributional kernels of analytic Fourier integral operators, cf.~\cite{MR4436039}. In fact, it is easy to see that we can follow \cite[Section 18.7]{MR4436039} in order to microlocalize the distributional kernels, i.e.~$\R^N$ is replaced by an open cone in Theorem \ref{Osc} 
with the obvious modification of $\Delta_\Phi$ etc. 
Then Theorem \ref{Osc} combined 
with Theorem \ref{WFkernel} gives the following statement.
\begin{Thm}\label{FIO-WF}
	Let $\Omega_1\subseteq\R^{n_1},\Omega_2\subseteq\R^{n_2}$ be two open sets
	and $\Gamma\subseteq\R^N\!\setminus\!\{0\}$ an open cone.
	If $\Phi$ is a real-valued real-analytic phase function
	on $\Omega_1\times\Omega_2\times\Gamma$ then any analytic Fourier integral operator $A$ 
	with phase $\Phi$
	satisfies 
	\begin{equation}
		\WF_{[\fM]}Au\subseteq\mathfrak{R}_\Phi(\WF_{[\fM]} u),
		\qquad \fa u\in \E^\prime(\Omega_2)
	\end{equation}
	for any [semiregular] weight matrix $\fM$, where
	\begin{equation*}
		\mathfrak{R}_\Phi(E)=\Set*{(x,\xi)\in T^\ast\Omega_1\setminus\{0\}\given \ex (y,\eta)\in E\;
			\text{with }
		(x,y,\xi,-\eta)\in V_\Phi(\Delta_\Phi)}
	\end{equation*}
	for any subset $E\subseteq T^\ast\Omega_2\!\setminus\{0\}$
	and $V_\Phi: \Delta_\Phi\rightarrow T^\ast\Omega_1\!\setminus\!\{0\}\times T^\ast\Omega_2\!\setminus\!\{0\}$ is the map given by
	$V_\Phi(x,y,\theta)=((x,y,d_x \Phi(x,y,\theta)),(x,y,-d_y\Phi(x,y,\theta)))$.
\end{Thm}

In the case of analytic pseudodifferential operators we give
some more details following mainly the presentation of 
\cite{MR4436039}: Let $\Omega\subseteq\R^n$ be an open set.
Recall that a pseudoanalytic amplitude $a$ of order $d\in\R$ on $\Omega$ is 
a smooth function $a\in\CC^\infty(\Omega\times\Omega\times\R^n)$ 
which can be extended to a smooth function
on $\Omega_{\C}\times\Omega_{\C}\times\R^n$, where $\Omega_{\C}$ is a neighborhood of $\Omega$ in $\C^n$, and $a$ is holomorphic with respect to $(z,w)\in\Omega_{\C}\times\Omega_{\C}$.
Moreover, for each compact subset $K\subseteq\Omega_{\C}\times\Omega_{\C}$ there are constants $C,\rho>0$ 
such that  for all $\gamma\in\Z_+^n$, all $(z,w)\in K$ and every $\xi\in\R^n$ with $\xit\geq\rho\max\{1,\abs{\gamma}\}$ we have that
\begin{equation}\label{Amplitude}
	\abs*{\partial^\gamma_{\xi}a(z,w,\xi)}\leq C^{\abs{\gamma}+1}\abs{\gamma}!\xit^{d-\abs{\gamma}}.
\end{equation}
We denote the space of pseudoanalytic amplitudes by $S_a(\Omega\times\Omega)$.
The analytic pseudodifferential operator $A:\E^\prime(\Omega)\rightarrow \D^\prime(\Omega)$ associated to the amplitude $a\in S_a(\Omega\times\Omega)$
is defined by
\begin{equation*}
	Au(x)=\frac{1}{(2\pi)^n}\int_{\R^n} e^{ix\xi}\left\langle u(y) , a(x,y,\xi)e^{-iy\xi}
	\right\rangle_y \,d\xi.
\end{equation*} 
The kernel of $A$ is the oscillatory integral
\begin{equation}\label{PseudoKernel}
	A(x,y)=\int_{\R^n}\! a(x,y,\xi)e^{i(x-y)\xi}\,d\xi.
\end{equation}
We observe that if we consider the phase-function $\Phi(x,y,\theta)=(x-y)\theta$
then we have that $\Delta_\Phi=\Set*{(x,y,\theta)\given x=y}$
and $d_{(x,y)}=\theta dx-\theta dy$.
Therefore Theorem \ref{WFkernel} combined with Theorem \ref{Osc}
implies that 
\begin{equation}\label{PseudoIncl}
	\WF_{[\fM]}Au\subseteq\WF_{[\fM]}u,\qquad \fa u\in\E^\prime(\Omega),
\end{equation}
for any analytic pseudodifferential operator $A$ and
any [semiregular] weight matrix $\fM$.

Following \cite[Chapter 17]{MR4436039} we can microlocalize the notion of analytic pseudodifferential operators:
If $\Omega\subseteq\R^n$ is an open set and
$\Gamma\subseteq\R^n\!\setminus\!\{0\}$ is an open cone then
a function $a\in\CC^\infty(\Omega\times\Omega\times\Gamma)$
is a pseudoanalytic amplitude if the following conditions hold: There is a complex neighborhood
$\Omega^\C$ of $\Omega$ such that for $\xi\in\Gamma$ fixed $a(x,y,\xi)$ extends to a holomorphic function on $\Omega^\C\times\Omega^\C$.
Each derivative $D^\alpha_x D^\beta_y D^\gamma_\xi a$ is bounded
near $\xi=0$ for any $(\alpha,\beta,\gamma)\in\Z_+^{3n}$
 and, most importantly, $a$ satisfies 
\eqref{Amplitude} for $(x,y,\xi)\in K\times K\times\Gamma$,
where $K$ is a compact subset of $\Omega^\C$.
We call a continuous operator $A:\E^\prime(\Omega)\rightarrow\D^\prime(\Omega)$ 
an analytic pseudodifferential operator associated to the open cone
$\Gamma\subseteq\R^n \setminus\{0\}$ if
there is a pseudoanalytic amplitude on $\Omega\times\Omega\times\Gamma$ such that
\begin{equation*}
	Au(x)=\frac{1}{(2\pi)^n}\int_\Gamma\!e^{ix\xi}\left\langle u(y) , a(x,y,\xi)e^{-iy\xi}
	\right\rangle_y \,d\xi.
\end{equation*}
Accordingly the distribution kernel of $A$ is given by
the oscillatory integral \eqref{PseudoKernel} with $\Gamma$ replacing $\R^n$ as the ``domain of integration''.
It is easily seen that the proof of Theorem  \ref{Osc} works also in this case and combined with Theorem \ref{WFkernel}
gives immediately the following statement, where $\Psi^d_a(\Omega;\Gamma)$ denotes the
space of analytic pseudodifferential operators of order $d$
associated to the cone $\Gamma$.
\begin{Thm}\label{MicroCHar}
	If $A\in\Psi^d_a(\Omega;\Gamma)$ then
	\begin{equation*}
		\WF_{[\fM]} Au\subseteq\WF_{[\fM]} u\cap\Omega\times\Gamma
	\end{equation*} 
	for all $u\in\E^\prime(\Omega)$ and all [semiregular]
	weight matrices $\fM$.
	
	In particular $A$ is $[\fM]$-pseudolocal, i.e. 
	$\singsupp_{[\fM]} Au\subseteq \singsupp_{[\fM]} u$
	for all $u\in\E^\prime(\Omega)$.
\end{Thm} 
Theorem \ref{MicroCHar} means in particular that, if we transition to sheaves of
analytic pseudodifferential operators and analytic microdifferential operators, then the sections of these sheaves
behave again well with respect to the ultradifferentiable setting. 
Before we give more details in this direction we want to 
give a more precise statement than Theorem \ref{MicroCHar}.
To this end, we need to introduce the weight function $\omega_\bM$ associated to a weight sequence $\bM$, which is given by $\omega_\bM(0)=0$ and
\begin{equation}\label{omegaM}
	\omega_\bM(t)=\sup_{k\in\Z_+}\log\frac{t^k}{M_k},\quad t>0.
\end{equation}
We notice that $\omega_\bM$ is a continuous function on
$[0,\infty)$ such that $\omega_\bM(t)\rightarrow\infty$ for 
$t\rightarrow\infty$.
In particular, $\omega(t)$ grows faster than $\log t^p$ for any
integer $p$, cf.~\cite{Komatsu1973}. %We can recover
%the weight sequence $\bM$ from its associated weight sequence:
%\begin{equation*}
%	M_k=\sup_{t\geq0}\frac{t^k}{\exp(\omega_\bM(t))}.
%\end{equation*}
Note also that 
\begin{equation}\label{ExpWeight}
	e^{-\omega_\bM(t)}=\inf_{k\in\Z_+}\frac{M_k}{t^k},\qquad t>0.
\end{equation}
\begin{Def}
	Let $\fM$ be a weight matrix and $a\in S_a(\Omega\times\Omega\times\Gamma)$.
	\begin{enumerate}
		\item The amplitude $a$ is $\{\fM\}$-regularizing in some
	open	conic set $V\times W\times\Gamma^\prime\subseteq \Omega\times\Omega\times\Gamma$ if there are 
	  a weight sequence $\bM\in\fM$
		and positive constants $C$, $h$ and $\rho$ such that
		\begin{equation}\label{M-Regularizing}
			\abs*{D^\alpha_x D^\beta_ya(x,y,\xi)}\leq Ch^{\abs{\alpha+\beta}}M_{\abs{\alpha+\beta}}\exp(-\omega_{\bM}(\rho\xit)),\qquad (x,y,\xi)\in V\times W\times\Gamma.
		\end{equation}
	\item $a$ is $(\fM)$-regularizing in $V\times W\times\Gamma^\prime\subseteq \Omega\times\Omega\times\Gamma$ if for every $\bM\in\fM$ and all $h,\rho>0$
	there is a constant $C>0$ such that \eqref{M-Regularizing} holds.
	\end{enumerate}
If $a$ is the amplitude associated to the pseudodifferential operator $A$ then we define the $[\fM]$-essential support $\essupp_{[\fM]} A$ of $A$ in the following way:
 $(x_0,\xi_0)\notin\essupp_{[\fM]} A$ if there is a conic neighborhood $V\times W\times\Gamma^\prime$ of $(x_0,x_0,\xi_0)$ such that $a$ is $[\fM]$-regularizing in $V\times W\times\Gamma$.
\end{Def}

\begin{Thm}\label{EssentialThm}
	Let $A\in\Psi_a(\Omega;\Gamma)$ and $\fM$ be a
	[semiregular] weight matrix. Then
	\begin{equation*}
		\WF_{[\fM]} Au\subseteq\WF_{[\fM]} u\cap\essupp_{[\fM]} A,\qquad \fa u\in\E^\prime(\Omega).
	\end{equation*} 
\end{Thm}
\begin{proof}
	In view of Theorem \ref{MicroCHar} it suffices to show
	that $(x_0,\xi_0)\notin\essupp_{[\fM]} A$ implies
	$(x_0,x_0,\xi_0,-\xi_0)\notin\WF_{[\fM]} A$.
	We return to the proof of Theorem \ref{Osc} with $x$ and $\theta$ replaced by $(x,y)$ and $\xi$, respectively, 
	and set $\Phi(x,y,\xi)=(x-y)\xi$. The proof of Theorem
	\ref{Osc} gives that given an open neigborhood $U$
	and a covering $F_j$, $j=1,\dotsc,\nu$
	 of open convex cones of $\Gamma$ we can write
	\begin{equation}\label{Repr}
		A\vert_U= h_0+\sum_{j=1}^\nu \bv_{F_j} G_j
	\end{equation}
	where $h_0$ is an analytic function in $U$
	and $G_j$ are holomorphic functions in $U+i\Gamma^j$
	(without loss of generality of slow growth) with
	$\Gamma^j=\Set*{(x^\prime,y^\prime)\in\R^{2n}\given (x^\prime-x^\prime)\xi>0\fa \xi\in F_j}$. In fact
	\begin{equation*}
		G(x+ix^\prime,y+iy^\prime)=\int e^{i(x+ix^\prime-y-iy^\prime)\xi}\chi(\xi)
		a(x+ix^\prime,y+iy^\prime,\xi)\,d\xi
	\end{equation*}
	where $0\leq\chi_j\leq 1$ and $\supp\chi\subseteq F_j\cap\Set{\xi\given \abs{\xi}\geq 1/2}$.
	If we write $h_j=\bv G_j\in\D^\prime(U)$ then
	obviously $\WF_{[\fM]}h_j\subseteq V\times (\Gamma^j)^\circ$. 
	We observe that (cf.~\cite[Subsection 7.2]{Liess1999})
	\begin{equation*}
		(\Gamma^j)^\circ=\overline{\Set*{(\xi,-\xi)\given \xi\in
			F_j}}.
	\end{equation*}
	Now let $(x_0,\xi_0)\notin\essupp_{[\fM]}A$. Then there exist a neighborhood $V\times W$ of $(x_0,x_0)$ and a open 
	cone $\Gamma^\prime$ containing $\xi_0$ such that
	$a(x,y,\xi)$ is $[\fM]$-regularizing in $V\times W\times\Gamma^\prime$; w.l.o.g.~ $V\times W
	\subseteq U$. We may also assume that there is some $j_0\in\Set{1,\dotsc,\nu}$ such that 
		$\xi_0\in F_{j_0}\subseteq\Gamma^\prime$
		and $\xi_0\notin F_j$ for $j\neq j_0$.
We have that
\begin{equation}\label{Integral1}
	h_{j_0}(x,y)=\int e^{i(x-y)\xi}\chi(\xi) a(x,y,\xi)\,d\xi
\end{equation}
as an oscillatory integral. However for $(x,y)\in V\times W$
the integral in the right-hand side of \eqref{Integral1}
is absolutely convergent since the integrand can be bounded by an
rapidly decreasing function. Because this bound is 
uniformly for $(x,y)\in V\times W$ we conclude that
$h_{j_0}\vert_{V\times W}\in\CC^\infty(V\times W)$.
We compute the derivatives of $h_{j_0}$:
\begin{equation*}
	\begin{split}
		D^\alpha_x D^\beta_y h_{j_0}(x,y)&=\int \chi(\xi) D^\alpha_xD^\beta_y
		\left(e^{i(x-y)\xi}a(x,y,\xi)\right)\,d\xi\\
		&=\int\chi(\xi)e^{i(x-y)\xi}\sum_{\substack{\alpha^\prime\leq\alpha\\
		\beta^\prime\leq \beta}}\binom{\alpha}{\alpha^\prime}
	\binom{\beta}{\beta^\prime} \xi^{\alpha^\prime}
	(-\xi)^{\beta^\prime} D^{\alpha-\alpha^\prime}_x
	D^{\beta-\beta^\prime}_y a(x,y,\xi)\,d\xi.
	\end{split}
\end{equation*}
If $A$ is $\{\fM\}$-regularizing 
we can therefore estimate using \eqref{ExpWeight} and \eqref{M-Regularizing} and obtain
\begin{equation*}
	\begin{split}
		\abs*{D^\alpha_x D^\beta_y h_{j_0}(x,y)}
		&\leq \int\limits_{\abs{\xi}\geq 1/2}
		\sum_{\substack{\alpha^\prime\leq\alpha\\
				\beta^\prime\leq \beta}}\binom{\alpha}{\alpha^\prime}
		\binom{\beta}{\beta^\prime} 
		\abs{\xi}^{\abs{\alpha^\prime+\beta^\prime}} 
		\abs{D^{\alpha-\alpha^\prime}_x
		D^{\beta-\beta^\prime}_y a(x,y,\xi)}\,d\xi\\
		&\leq C\int\limits_{\abs{\xi}\geq 1/2}
		\sum_{\substack{\alpha^\prime\leq\alpha\\
			\beta^\prime\leq\beta}}
		\binom{\alpha}{\alpha^\prime}\binom{\beta}{\beta^\prime}
		h^{\abs{\alpha+\beta}-\abs{\alpha^\prime+\beta^\prime}}
		M_{\abs{\alpha+\beta}-\abs{\alpha^\prime+\beta^\prime}} 
		\abs{\xi}^{\abs{\alpha^\prime+\beta^\prime}}
		e^{(-\omega_\bM(\rho\abs{\xi}))}\\
		&\leq C\int\limits_{\abs{\xi}\geq 1/2}
		\sum_{\substack{\alpha^\prime\leq\alpha\\
				\beta^\prime\leq\beta}}
		\binom{\alpha}{\alpha^\prime}\binom{\beta}{\beta^\prime}
		h^{\abs{\alpha+\beta}-\abs{\alpha^\prime+\beta^\prime}}
		M_{\abs{\alpha+\beta}-\abs{\alpha^\prime+\beta^\prime}}\times\\
		&\qquad\qquad\qquad\times
		M_{\abs{\alpha^\prime+\beta^\prime}+n+1}
		\rho^{\abs{\alpha^\prime+\beta^\prime}+n+1}
		\int\limits_{\abs{\xi}\geq 1/2}\abs{\xi}^{-n-1}\,d\xi\\
		&\leq C \rho^{n+1} \left(2\max\{h,\rho\}\right)^{\abs{\alpha+\beta}}
		M_{\abs{\alpha+\beta}+n+1}
		\int_{\abs{\xi}\geq 1/2}\abs{\xi}\,d\xi.
	\end{split}
\end{equation*}
Using \eqref{R-DerivClosed} $n+1$-times we infer that
there are a weight sequence $\bM^\prime\in\fM$,
constants $C^\prime,h^\prime>0$ such that
\begin{equation*}
	\sup_{(x,y)\in V\times W}\abs{D^\alpha_xD^\beta_y h_{j_0}(x,y)}\leq C^\prime (h^\prime)^{\abs{\alpha+\beta}}M^\prime_{\abs{\alpha+\beta}},
	\qquad \fa (\alpha,\beta)\in\Z_+^{2n}.
\end{equation*}
That means in particular that $h_{j_0}\vert_{V\times W}\in
\Rou{\fM}{V\times W}$ and therefore using \eqref{Repr}
we conclude that $(x_0,x_0,\xi_0,-\xi_0)\notin\WF_{\{\fM\}} A$.
	
	In the Beurling case, it is easy to see how to
	adapt the above proof to show the implication
	\begin{equation*}
		(x_0,\xi_0)\notin\essupp_{(\fM)}\Longrightarrow
		(x_0,x_0,\xi_0,-\xi_0)\notin \WF_{(\fM)}A
	\end{equation*} 
	if $\fM$ is $B$-semiregular.
\end{proof}

We have seen that the action of analytic differential operators defined on analytic manifolds
on semiregular ultradifferentiable classes is well behaved. 
Theorem \ref{WFkernel} will show that the same is true for analytic pseudodifferential operators.
It will be convenient to first introduce the notion of microdistributions of ultradifferentiable degree.
For the definition of microdistributions in the smooth and analytic category we refer to
\cite{MR0597144} and \cite{MR4436039}.
The extension to the ultradifferentiable category is straight-forward:

%If $\Omega\subseteq\R^n$ is an open set, $\fM$  a [semiregular] weight matrix 
%and $(x_0,\xi_0)\in\Omega\times\R^n\setminus\{0\}$ then we define on $\D^\prime(\Omega)$
%an equivalence relation:
If $\fM$ is a  weight matrix and
if $\mathcal{U}\subseteq\R^n\times\R^n\!\setminus\!\{0\}$ is a conic open set then
 we define on $\D^\prime(\pi_1(\mathcal{U}))$
an equivalence relation:
\begin{equation*}
	%u\sim_{[\fM]} v\text{ at }(x_0,\xi_0)\quad\Longleftrightarrow\quad (x_0,\xi_0)\notin\WF_{[\fM]}(u-v).\\	
	u\sim_{[\fM]}v\text{ in }\mathcal{U}\quad \Longleftrightarrow\quad
	\WF_{[\fM]}(u-v)\cap \mathcal{U}=\emptyset.
\end{equation*}
For open conic sets $\mathcal{V}\subseteq\mathcal{U}\subseteq T^\ast\R^n\!\setminus\{0\}$ there is a well-defined map 
\begin{equation*}
\rho_{\mathcal{V},\mathcal{U}}:\quad\D^\prime(\pi(\mathcal{U}))/\!\sim_{[\fM]}^{\mathcal{U}}\longrightarrow \D^\prime(\pi(\mathcal{V}))/\!\sim_{[\fM]}^{\mathcal{V}}
\end{equation*}
induced by the restriction map $\D^\prime(\pi(\mathcal{U}))\rightarrow\D^\prime(\pi(\mathcal{V}))$.
Thence $(\D^\prime(\pi(\mathcal{U}))/\!\sim_{[\fM]}^{\mathcal{U}},\rho_{\mathcal{V},\mathcal{U}})$
is a presheaf on $T^\ast\R^n\!\setminus\!\{0\}$. We denote the resulting sheaf by $\mathsf{D}^{micro}_{[\fM]}(\R^n)$ and call it the sheaf of microdistributions of degree $[\fM]$.
 The stalk of this sheaf at a point $(x,\xi)\in T^\ast(\R^n)$
is equivalent to $\D^\prime_{x_0}/\sim^{(x_0,\xi_0)}_{[\fM]}$
	where $\D^\prime_{x_0}$ is the set of germs of distributions at $x_0$ and, if $u,v\in\D^\prime_{x_0}$, then
%	$u\sim_{(x_0,\xi_0)v$
\begin{equation*}
u	\sim_{[\fM]}^{(x_0,\xi_0)} v\quad\Longleftrightarrow\quad (x_0,\xi_0)\notin\WF_{[\fM]}(u-v).
\end{equation*}
If we denote the sheaf of (smooth) microdistributions
and the sheaf of analytic microdistributions by
$\mathsf{D}^{micro}(\R^n)$ and $\mathsf{D}^{micro}_a(\R^n)$,
respectively, then, by Proposition \ref{WFproperties}, there are sheaf morphisms
\begin{equation*}
	\mathsf{D}_a^{micro}(\R^n)\longrightarrow
	\mathsf{D}^{micro}_{(\fM)}(\R^n)\longrightarrow
	 \mathsf{D}^{micro}_{\{\fM\}}(\R^n)\longrightarrow\mathsf{D}^{micro}(\R^n)
\end{equation*} 
for any weight matrix $\fM$ satisfying \eqref{AnalIncl}.
Now assume that $\fM$ and $\fN$ are two weight matrices.
Then we set
\begin{align*}
	\fN &\{\preceq\}\fM & &:\Longleftrightarrow
	& \fa \bN\in\fN\;\ex\bM\in\fM\;\ex h>0\;\ex C>0\quad
	N_k&\leq Ch^kM_k\qquad \fa k\in\Z_+,\\
	\fN&(\preceq)\fM & &:\Longleftrightarrow 
	&\fa \bM\in \fM \;\ex \bN\in\fN\; \ex h>0\;\ex C>0\quad
	N_k&\leq Ch^k M_k\qquad \fa k\in\Z_+,\\
	\fN &\{\lhd)\fM & &:\Longleftrightarrow &
	\fa \bN\in\fN\;\fa\bM\in\fM\;\fa h>0\;\ex C>0\quad
	N_k &\leq Ch^kM_k\qquad \fa k\in\Z_+.
\end{align*}
If $\fN[\preceq]\fM$ then $\WF_{[\fM]}u\subseteq\WF_{[\fN]}u$ for all $u\in\D^\prime$. Moreover,
$\WF_{(\fM)}u\subseteq\WF_{\{\fN\}}u$ if $\fN\{\lhd)\fM$.
Therefore there are sheaf morphisms
$\mathsf{D}^{micro}_{[\fN]}(\R^n)\rightarrow \mathsf{D}^{micro}_{[\fM]}(\R^n)$, when $\fN[\preceq]\fM$,
and $\mathsf{D}^{micro}_{\{\fN\}}(\R^n)\rightarrow \mathsf{D}^{micro}_{(\fM)}(\R^n)$ for $\fN\{\lhd)\fM$.

We denote the space of (continuous) sections of $\mathsf{D}_{[\fM]}^{micro}(\R^n)$ over an open conic
set $\mathcal{U}\subseteq T^\ast \R^n$ by $\mathfrak{D}^{micro}_{[\fM]}(\cU)$.
Note that we have an obvious notion of ``support'' for elements of $\mathfrak{D}^{micro}_{[\fM]}(\cU)$.
Namely, if $(x,\xi)\in\cU$ then $(x,\xi)\notin \supp_{[\fM]} [u]$, $[u]\in\mathfrak{D}^{micro}_{[\fM]}(\cU)$,
if there is a germ $v$ of a distribution at $x_0$ such that $v$ is a representative of $[u]$ at $(x,\xi)$ and $v\sim_{[\fM]} 0$ at $(x,\xi)$. 
We call an operator $T: \mathfrak{D}^{micro}_{[\fN]}(\cU)\rightarrow\mathfrak{D}^{micro}_{[\fM]}(\cU)$ local if
\begin{equation*}
	\supp_{[\fM]} T[u]\subseteq \supp_{[\fN]}[u]\qquad \fa[u]\in
	\mathfrak{D}^{micro}_{[\fN]}(\cU).
\end{equation*}
Obviously we can also introduce similar notions for
operators $S:\,\mathfrak{D}^{micro}_{(\fN)}\rightarrow\mathfrak{D}^{micro}_{\{\fM\}}$, $R:\,\mathfrak{D}^{micro}_{\{\fN\}}\rightarrow
\mathfrak{D}^{micro}_{(\fM)}$, etc.
It follows that all the operators induced by the sheaf morphisms
above are local.

We turn our attention to the sheaf of microdifferential operators
on $T^\ast\R^n$. For the definition and properties see \cite[Section 17.5]{MR4436039}.
%We note that 
%
%It follows that $A$ acts naturally on sections of $\mathsf{D}^{micro}_{[\fM]}(\R^n)$.
Following the arguments in \cite[pp.~278-285]{MR0597145} (cf.~ also \cite{MR4436039})
we see that sections of the sheaf $\Psi_a^{micro}$ of analytic microdifferential operators
are acting on sections of $\mathsf{D}^{micro}_{[\fM]}(\R^n)$ ($\fM$ is always a [semiregular]
weight matrix) and they do so in a "local" manner:
Let $\mathcal{U}$ be a conic subset of $T^\ast\R^n$, $\fM$ be a [semiregular] weight matrix,
$[A]\in\Psi^{micro}_a(\mathcal{U})$ and $[u]\in\mathfrak{D}^{micro}_{[\fM]}(\mathcal{U})$;
then	$\supp [A][u]\subseteq [u]$
by Theorem \ref{MicroCHar}. In fact Theorem \ref{EssentialThm}
gives a more precise statement.
We note that we can easily extend the notion of $\essupp_{[\fM]}$
to analytic microdifferential operators $[A]$: $(x,\xi)\notin\essupp_{[\fM]}[A]$ if there is a representative
$A$ of $[A]$ defined in a conic neighborhood of $(x,\xi)$
such that $(x,\xi)\notin \essupp_{[\fM]}[A]$.
\begin{Thm}
	If $[A]\in\Psi^{micro}_a(\mathcal{U})$ then
	\begin{equation*}
		\supp_{[\fM]} [A][u]\subseteq \supp_{[\fM]} [u]\cap \essupp_{[\fM]}[A]\qquad \fa [u]\in\mathfrak{D}^{micro}_{[\fM]}(\mathcal{U})
	\end{equation*}
	for all [semiregular] weight matrices $\fM$.
\end{Thm}

Similarly we can extend the notion of ellipticity to microdifferential operators. Since any elliptic analytic pseudodifferential operator has locally in the algebra of analytic pseudodifferential operators an inverse modulo analytic regularizing operators (\cite[Proposition 17.2.21]{MR4436039}) we obtain 
\begin{Thm}\label{EllipticRegThm}
	Let $A\in\Psi_a^{micro}(\mathcal{U})$ be elliptic.
	Then
	\begin{equation*}
		\supp_{[\fM]} [A][u]=\supp_{[\fM]}[u]\qquad \fa [u]\in\mathfrak{D}^{micro}_{[\fM]}(\mathcal{U})
	\end{equation*}
	for every [semiregular] weight matrix $\fM$.
\end{Thm}
This proves Theorem \ref{StandardElliptRegThm}.
We finish this section by turning our attention to classic analytic pseudodifferential operators. 
For the definition and basic theorey we refer to
\cite[Chapter 17]{MR4436039} and note here
that the symbol $a$ associated to  a classical analytic pseudodifferential operator $A$ of order $d$ on $\Omega\times\Gamma$ is of the form
\begin{equation*}
	a(x,\xi)\sim a_d(x,\xi)+a_{d-1}(x,\xi)+a_{d-2}(x,\xi)+\dots
\end{equation*}
where the symbols $a_{d-j}(x,\xi)$ are homogeneous of degree $d-j$ in the second variable.
We call $a_d$ the principal symbol of $A$ and
can therefore define the characteristic set of $A$:
\begin{equation*}
	\Char A=\Set*{(x,\xi)\in T^\ast\Omega\setminus\{0\}\given a_{d}(x,\xi)=0}.
\end{equation*}
Hence $A$ is elliptic in the set $\mathcal{U}\setminus\Char A$, cf.~\cite[Section 17.5]{MR4436039}.
We say that a section $[A]\in\Psi_a^{micro}(\mathcal{U})$ is a classical microdifferential operator
on $\mathcal{U}$ if there is a representative $A$ of $[A]$ 
which is a classical pseudodifferential operator in $\mathcal{U}$.
By \cite[17.5.2]{MR4436039} the inverse of an elliptic classical microdifferential operator $[A]\in\Psi^{micro}_a(\mathcal{U})$ is again a classical microdifferential operator.
The following result is a direct consequence of Theorem \ref{EllipticRegThm}.
\begin{Thm}
	Let $A$ be a classical analytic pseudodifferential operator on some open set $\Omega\subseteq\R^n$ and $\fM$ be a [semiregular] weight matrix.
	Then
	\begin{equation*}
		\WF_{[\fM]} u=\WF_{[\fM]}Au\cup\Char A\qquad\fa u\in\E^\prime(\Omega).
	\end{equation*}
\end{Thm}

Note that all these considerations hold also if we replace $\R^n$
by a paracompact real-analytic manifold $\M$, cf.~\cite[Chapter 17]{MR4436039}.

\section{Hypoelliptic differential operators of principal type} 
\label{Sec:principal}
In order to formulate the main results of the last two sections we need to recall
some facts and notations from the modern geometric theory of linear partial differential equations
(with real-analytic coefficients).

%\begin{Def}
%	Let $A$ be a classical pseudodifferential operator of order $d$ in $\Omega$ 
%	with principal symbol $a_d$.
%	We say that $A$ is of principal type if 
%	\begin{equation*}
%%		d a_d\wedge\sigma (x_0,\xi_0)\neq 0
%	\end{equation*}
%for all $(x_0,\xi_0)\in\Char A$. Here $\sigma=\xi dx=\sum \xi_jdx_j$ denotes the tautological one-form. 
%\end{Def}
If $f\in\CC^{\infty}(T^\ast\Omega\!\setminus\!\{0\})$ is a real-valued smooth function, %homogeneous in the second variable and real-valued 
then the Hamiltonian vector field associated to $f$
is 
\begin{equation*}
	H_f= \sum_{j=1}^n\left(\frac{\partial f}{\partial \xi_j}\frac{\partial}{\partial x_j}
	-\frac{\partial f}{\partial x_j}\frac{\partial}{\partial \xi_j}\right).
\end{equation*}
Now if $A$ is a classical pseudodifferential operator in $\Omega$ with principal symbol 
$a_d=b+ic$ where $b$ and $c$ are real-valued, then a null bicharacteristic leaf of $A$ is
an integral manifold of the Lie algebra generated by $b$ and $c$, which is contained in $\Char A$. 
If the dimension of a null bicharacteristic leaf is $1$ then we refer to it as a null bicharacteristic curve.
\begin{Def}[{\cite[Definition 23.1.37]{MR4436039}}]
	We say that a differential operator $P$ with
	principal symbol $p_d$ satisfies Condition $(\mathbf{P})$ if
	$\real p_d$ is linear independent of the radial vector field
	on $T^\ast\Omega\setminus\{0\}$ on each point 
	$x\in\{y\in T^\ast\Omega\setminus\{0\}\real p_d(y)=0\}$
	and the restriction of $\imag p_d$ to any null bicharacteristic curve of $\real p_d$ does not change sign.
\end{Def}
We note here that a smooth hypoelliptic differential operator of principal type satisfies Condition $(\mathbf{P})$
in $T^\ast\Omega\!\setminus\!\{0\}$, cf.~\cite[Proposition 24.2.2]{MR4436039}. On the other hand
Treves \cite{MR0296509}, see also \cite[Theorem 24.2.1]{MR4436039}, proved that an analytic 
differential operator $P$ of principal type is smooth hypoelliptic or analytic hypoelliptic if and only if the following condition is satisfied:
\begin{equation}\tag{\textbf{Q}}\label{ConditionQ}
\text{There are no null bicharacteristic leaves of }P\text{ of positive dimension in }T^\ast\Omega.
\end{equation}
Thence Theorem \ref{PrincipalThm1} is a consequence of the following statement.
\begin{Thm}\label{GeneralTrevesThm}
	Let $P$ be an analytic differential operator of principal type
	which satifies Condition $(\mathbf{P})$. 
	For any [semiregular] weight matrix $\fM$
 the following statements are equivalent:
	\begin{enumerate}
	\item $P$ satisfies Condition \eqref{ConditionQ}.
	\item $P$ is $[\fM]$-hypoelliptic.
	\end{enumerate}
\end{Thm}
 \begin{proof}
 	Suppose \eqref{ConditionQ} holds and let $x_0\in\pi(\Char P)$.
 	By \cite[(24.6.31)]{MR4436039} we have in a suitably small neighborhood $U$ of $x_0$
 	\begin{equation}\label{ParametrixPrincipal}
 		P(x,D_x)Au=Eu+Fu
 	\end{equation}
 for all $u\in\CC^\infty_0(U)$ where the distributional kernel of the operator $A:\E^\prime\rightarrow\D^\prime$ 
 is analytic off the diagonal,
 $E$ is an elliptic analytic pseudodifferential operator and $F$ is analytic regularizing.
 %Since the operator $A$ is also semiregular in the sense of \cite{MR4436039}*{Definition 2.3.1} we can extend 
 %the formula 
 
 Let $U_0\Subset U$ and $\chi\in\CC^\infty_0(U)$ with $\chi\vert_{U_0}\equiv 1$.
 Then
 \begin{equation*}
 	P(x,D)A \chi E^{-1}u=E \chi E^{-1} u +F_1u
 \end{equation*}
 where $F_1=F\chi E^{-1}$.
Thus
\begin{equation*}
	\left(P A \chi E^{-1} u\right)\vert_{U_0}= u+ F_2u
\end{equation*}
 for $u\in\CC^\infty_0(U_0)$, cf.~\cite[Proposition 17.2.21]{MR4436039}.
 Here $F_1,F_2$ are analytic regularizing in $U_0$, i.e.~$F_j(\E^\prime(U_0))\vert_{U_0}\subseteq\An(U_0)$.
 We put 
 \begin{equation*}
 	G=\left(A\chi E^{-1}\right)\vert_{U_0}.
\end{equation*}
 By transposition  we obtain
 \begin{equation*}
 	G^{\top}P^{\top}u=u +F_3u
 \end{equation*}
 for all $u\in\E^\prime(U_0)$.
 where $F_3$ is again analytic regularizing.
Now note that the distribution kernels of $A^{\top}$ and
$(E^{-1})^{\top}$ are also analytic off the diagonal.
By Theorem \ref{WFkernel} it therefore follows that $G^{\top}$ is $[\fM]$-pseudolocal:
\begin{equation*}
	\singsupp_{[\fM]} G^{\top}u\subseteq\singsupp_{[\fM]} u
\end{equation*}
for all $u\in\E^\prime(U_0)$.
This implies  that $\singsupp_{[\fM]} P^{\top}u=\singsupp_{[\fM]} u$ for all $u\in\E^\prime(U_0)$.
Thence  $P^{\top}$ is $[\fM]$-hypoelliptic at $x_0$, i.e.
\begin{equation*}
	x_0\notin\singsupp_{[\fM]}P^{\top}u\;\Longrightarrow\; x_0\notin\singsupp_{[\fM]}u
\end{equation*}
for all distributions $u$ defined near $x_0$.
Since \eqref{ConditionQ} is invariant under transposition we have proven the implication 
$(1)\Rightarrow (2)$.

On the other hand suppose that \eqref{ConditionQ} is violated at some point $x_0$.
If we argue in a small enough neighborhood $U$ of $x_0$ then
we can  take the phase function $\varphi$ and the formal solution
$a(z,\lambda)=\sum_{j=0}^\infty e^{i\lambda\varphi(z)}\lambda^{-j}a_j(z,\lambda)$ of $P(z,D_z)a(z,\lambda)$ from \cite[Sections 24.3 \& 24.5]{MR4436039}, see also \cite[Chapter 23]{MR4436039}. Here we extend the coefficients of $P$ holomorphically to a complex neighborhood
$\widetilde{U}$ of $U$ in $\C^n$.
In particular $\varphi\vert_{z_n=0}=z^\prime \xi_0^\prime+ i\lvert z^\prime\rvert^2$.
Moreover, for $N>0$, consider the approximate solution of the equation $Pu=0$ given by 
$e^{i\lambda\varphi}a^{[N]}$ where
\begin{equation*}
	a^{[N]}(z,\lambda)=\sum_{j\leq\frac{\lambda}{N}}\lambda^{-j}a_j(z,\lambda).
\end{equation*}

By \cite[Lemma 24.5.1 \& (24.5.7)]{MR4436039} we have, for $U$ small enough, that for
any compact set $K\subseteq U$ there is a constant $C_K>0$ such that 
\begin{equation}\label{Estimate1}
\fa\lambda\geq 1:\quad	\sup_{x\in K} \abs*{e^{i\lambda\varphi}
a^{[N]}(x,\lambda)}\leq C_K
\end{equation}
for $N\geq C_K$.
Moreover, \cite[Proposition 24.5.2 \& (24.5.)]{MR4436039} indicates that
 for any large $N$ there are
an open neighborhood $U^\C_N\subseteq\C^n$ of the origin
 and $C_N>0$ such that
 \begin{equation}\label{Estimate2}
 	\lambda\geq 1:\quad \max_{K_1}\;\abs*{P(z,D_z)e^{i\lambda\varphi a^{[N]}}}\leq C_N
\end{equation}
for any compact set $K_1\subseteq U^\C$. 
%Note that $U^\C_{N+1}\subseteq U^\C_N$.

%We fix now $N$, $U$ and $U^\C$ such that both \eqref{Estimate1} and \eqref{Estimate2} hold.

We recall from \cite[pp.~995--996]{MR4436039}
that $\partial_{x_j}\varphi(0)=0$ for $j=1,\dotsc,n$
and by the initial conditions $\partial^{\alpha}\varphi(0)=1$
for $\alpha\in\Z_+^n$ with $\abs{\alpha}=0$ and $\alpha_n=0$.

We observe for any $C^2$-functions $f,g$ we have the following formula:
\begin{equation*}
	\begin{split}
	D_jD_k\bigl(e^{i\lambda f}g\bigr)&=
	-i\lambda\bigl(\partial_j\partial_k f\bigr)e^{i\lambda f}g
	+ \lambda^2\bigl(\partial_j f\bigr)(\partial_k f)
	e^{i\lambda f}g\\
	+\lambda (\partial_k f)e^{i\lambda f}D_jg
	+\lambda(\partial_j f)e^{i\lambda f}D_kg
	+e^{i\lambda}D_jD_kg.
	\end{split}
\end{equation*}
We set now $f=\varphi$ and $g=a^{[N]}$ for abritrary $N$ and
obtain that
\begin{equation*}
	D^\alpha \bigl(e^{i\lambda\varphi(0)}a^{[N]}(0)\bigr)=-i\lambda \partial^\alpha\varphi(0) +D^\alpha a^{[N]}(0)
\end{equation*}
for any $\alpha\in\Z_+^n$ with $\abs{\alpha}=2$.
If now $K\subseteq U$ is any compact set containing the origin we
conclude that
\begin{equation}\label{LowerBound2}
	\begin{split}
	\sum_{\abs{\alpha}=2}\sup_{x\in K}\abs*{D^\alpha\bigl(e^{i\lambda\varphi(x)}a^{[N]}(x)\bigr)}
	&\geq \sum_{\abs{\alpha}=2}
	\abs*{\lambda \partial^\alpha\varphi(0)+D^\alpha a^{[N]}(0)}\\
	&\geq
	\sum_{\abs{\alpha}=2}\left(\lambda\abs*{\partial^\alpha\varphi(0)}
	-\abs*{D^\alpha a^{[N]}(0)}\right)\geq c_1\lambda-c_2
\end{split}
\end{equation} 
for all $\lambda\geq 1$
and some constants $c_1,c_2>0$.

We fix now $N$, $U$ and $U^\C$ such that both \eqref{Estimate1} and \eqref{Estimate2} hold. 
If we assume now that $P$ is $[\fM]$-hypoelliptic in $U$
than in particular the assumptions of Proposition \ref{NecProp}
are satisfied and therefore inequality \eqref{HypoEstimate} below holds
for $\Omega=U$ and $\widetilde{{\Omega}}=U^\C$.
If we set $u_\lambda=e^{i\lambda\varphi a^{[N]}}$ then
the right-hand side of \eqref{HypoEstimate} is bounded when $\lambda\rightarrow \infty$ for any choice of $K^\prime$ and $\widetilde{K}$.
But by \eqref{LowerBound2} we have that the left-hand side
of \eqref{HypoEstimate} is unbounded when $\lambda\rightarrow\infty$ for any $K$ containing the origin.
Hence \eqref{HypoEstimate} is violated and therefore $P$ cannot be $[\fM]$-hypoelliptic for any [semiregular] weight matrix $\fM$.
 \end{proof}
 \begin{Prop}\label{NecProp}
 	Let $P$ be a differential operator with analytic coefficients in a connected set 
 	$\Omega\subseteq\R^n$ and
 	assume that all coefficients of $P$ extends to the open 
 	set $\widetilde{\Omega}\subseteq\C^n$.
 	 We	suppose that for any $u\in\CC(\Omega)$ with $Pu\in\CC^\omega(\Omega)$
 	we have that $u\in\CC^2(\Omega)$.
 	Then for all compact set $K\subseteq \Omega$ there are compact sets $K^\prime\subseteq\Omega$, $\widetilde{K}\subseteq\widetilde{{\Omega}}$ and a constant
 	$C>0$ such that
 	\begin{equation}\label{HypoEstimate}
 		\sum_{\alp\leq 2}\sup_{x\in K}\abs*{D^\alpha u(x)}\leq 
 		C\left\{\sup_{x\in K^\prime}\abs*{u(x)}+\sup_{z\in \widetilde{K}}\abs*{P(z,D)u(z)}\right\}
 	\end{equation}
 	for all $u\in\mathcal{O}(\widetilde{\Omega})$. 
 \end{Prop}
 \begin{proof}
 	Following \cite[p.~37]{MR4436039} we confer $\CC^\omega(\Omega)$ with the following topology:
 	Let $(\Omega_\nu)_{\nu}$ be a sequence of neighborhoods of $\Omega$ in $\C^n$ such that
 	$\Omega_{\nu+1}\subseteq\Omega_{\nu}$ and $\bigcap_{\nu=1}^\infty\Omega_\nu=\Omega$.
 	Moreover, we may assume that all $\Omega_\nu$ are connected Runge open sets and that $\Omega_1\subseteq\widetilde{\Omega}$.
 	Obviously $\mathcal{O}(\Omega_\nu)\subseteq\mathcal{O}(\Omega_{\nu+1})$
 	and there is a canonical injection $\lambda_\nu:\mathcal{O}(\Omega_\nu)\rightarrow \CC^\omega(\Omega)$ for each $\nu\in\N$. In fact,
 	\begin{equation*}
 		\CC^\omega(\Omega)=\bigcup_{\nu=1}^\infty \lambda_\nu\left(\mathcal{O}(\Omega_\nu)\right)
 	\end{equation*}
 	as sets. We equip $\CC^\omega(\Omega)$ with the inductive limit topology coming from
 	the spectrum $(\mathcal{O}(\Omega_\nu),\lambda_\nu)$. We may denote this topological space
 	by $\An(\Omega)$. In particular, $\An(\Omega)$ is a  webbed space by \cite{MR551623} and a seminorm
 	$p$ is continuous on $\An(\Omega)$ if and only if $p\circ \lambda_\nu$ is a continuous seminorm on 
 	$\mathcal{O}(\Omega_\nu)$ for all $\nu\in\N$ according to \cite{MR226355}.
 	
 	It follows instantly that the space $E=\CC(\Omega)\times\An(\Omega)$ is also a webbed space, see \cite{MR551623}.
 	We claim  that the subspace $F=\Set{(u,f)\in E\given Pu=f}$ is webbed, too.
 	For this it is enough to show that $F$ is a sequentially closed subspace of $E$ by \cite{MR551623}: 
 	Let $(u_j,f_j)\in F$ be a sequence which converges in $E$, i.e.\
 	$(u_j,f_j)\rightarrow (u,f)$, then $u_j\rightarrow u$ in $\D^\prime(\Omega)$, which immediately gives that $Pu_j\rightarrow Pu$ in $\D^\prime(\Omega)$.
 	On the other hand we have also that $Pu_j=f_j\rightarrow f$ in $\D^\prime(\Omega)$ and
 	therefore $Pu=f$ in the sense of distributions, which means that $(u,f)\in F$.
 	
 	By assumption the map $\lambda: F\rightarrow \CC^2(\Omega)$ given
 	by $\lambda(u,f)=u$ is well-defined.
 	Moreover, note that if $(u_j,f_j)\rightarrow 0$ in $F$ and $\lambda(u_j,f_j)\rightarrow w$
 	in $\CC^1(U)$ then $w=0$. Thus $\lambda$ has closed graph.
 	Now the closed-graph theorem in the form given in \cite{MR551623} applies and therefore $\lambda$ is continuous. 
 \end{proof}
 
 \begin{Rem}
We observe that the proof of Theorem \ref{GeneralTrevesThm} implies that if $P$ is an analytic differential operator which satisfies (\textbf{P}) but not \eqref{ConditionQ}
then near each point $x_0\in\Omega$ there is a continuous function $u$
defined in a neighborhood $U$ of $x_0$ such
that $Pu$ is real-analytic in $U$, but $u$ is not $\mathcal{C}^2$
near $x_0$.
 \end{Rem}

\section{Metivier Operators}\label{Sec:Metivier}
In this section we show how the pattern of proof in \cite{MR0597752} can be modified to prove
Theorem \ref{MetivierThm}. Note that the proof on the operator side is the same as in the analytic category, only involving analytic pseudodifferential and analytic 
Fourier integral operators. 
We need only the that these operators do transform ultradifferentiable wavefronts sets in the correct way, which 
is true by Theorem \ref{FIO-WF} and Theorem \ref{MicroCHar},
respectively.

In fact, will work in the setting of Okaji \cite{MR0807494} which generalizes the setting of Metivier.
We say that an analytic classical pseudodifferential operator $P$  of order $d$ in an open set
$\Omega\subseteq\R^n$ with total symbol
\begin{equation*}
	p(x,\xi)=p_d(x,\xi)+p_{d-1}(x,\xi)+\dots
\end{equation*}
defined in $\Omega\times\Gamma$, where $\Gamma\subseteq\R^n\!\setminus\!\{0\}$ is an open cone
is an Metivier-Okaji Operator if the following conditions are satisfied:
\begin{enumerate}
	\item There are two analytic conic manifolds $\Sigma_1,\Sigma_2\subseteq T^\ast\Omega\!\setminus\!\{0\}$
	with the same codimension $\nu$, which are regular involutive, $\Sigma=\Sigma_1\cap \Sigma_2$ is
	a conic, symplectic,  analytic submanifold of $T^\ast\Omega\setminus\{0\}$ of codimension $2\nu$ and for each $\rho\in\Sigma$ we have that $T_\rho\Sigma=T_\rho \Sigma_1\cap T_\rho\Sigma_2$.
	\item For each point $\rho=(x_0,\xi_0)\in \Sigma$ there exists a conic neighborhood $\cU\subseteq T^\ast\Omega\setminus\{0\}$ of $\rho$ and $M,\mu\in\N$ such that 
	\begin{align*}
		\frac{\abs*{p_{d-j}(x,\xi)}}{\xit^{d-j}}
		&\leq C\left(d_{\Sigma_1}(x,\xi)+\left(d_{\Sigma_2}(x,\xi)\right)^{\mu}\right)^{M-j(\mu+1)/\mu}
		\\
		\frac{\abs*{p_d(x,\xi)}}{\xit^d}
		&\geq C^{-1}\left(d_{\Sigma_1}(x,\xi)+\left(d_{\Sigma_2}(x,\xi)\right)^\mu\right)^{M}
	\end{align*}
	for $(x,\xi)\in\cU\cap\Set{(x,\xi)\given\xit\geq 1}$ with $C$ being a constant only depending on $\cU$.
	Here $d_{\Sigma_j}(x,\xi)$, $(j=1,2)$, denotes a distance function between
	$(x,\xi/\xit)$ and $\Sigma_j\cap\{\xit=1\}$.
	\item The operator $P$ is hypoelliptic in $\Omega$ with loss of $M\mu/(\mu+1)$ derivatives.
\end{enumerate}
Thus Theorem \ref{MetivierThm} is a special case of the following theorem.
\begin{Thm}\label{MainMetivier}
	Let $P$ be a Metiver-Okaji operator in $\Omega$. Then $P$ is $[\fM]$-hypoelliptic for any [semiregular] weight matrix $\fM$.
\end{Thm}
If we follow the arguments in \cite[pp.~490--491]{MR0807494} then
we see that, using an elliptic analytic Fourier integral operator $F$, we can assume that
$x_0=0$, $\xi_0=e_n$,  such that 
the transformed operator $\tilde{P}=FPF^{-1}$ has the form 
\begin{equation}\label{canonicalform1}
	\tilde{P}=\sum_j\sum_{(\alp/\mu)+\bet=M-j(\mu+1)/\mu}c_{\alpha\beta}{x^\prime}^\alpha D^\beta_{x^\prime}
\end{equation}
where $x^\prime=(x_1,\dotsc,x_\nu,0,\dotsc,0)$ and $\xi_\nu=(\xi_1,\dotsc,\xi_\nu,0,\dotsc,0)$, ($\nu\leq n$) and the $c_\alpha\beta$ are classical analytic pseudodifferential operators
of order $d+(\alp-\bet-\mu M)/(\mu+1)$. The sets $\Sigma_j$, $j=1,2$, are transformed to
$\Sigma_1=\Set{(x,\xi)\in T^\ast\Omega\setminus\{0\}\given \xi^\prime=0}$, 
$\Sigma_2=\Set{(x,\xi)\in T^\ast\Omega\setminus\{0\}\given x^\prime=0}$.

Our assumptions on $P$ lead to the following conditions:
\begin{equation}\label{TransCond1}
	\abs*{y^\prime}+\abs*{\eta^\prime}\neq 0 \quad\Longrightarrow\quad	
	\sum_{\alp/\mu+\bet=M} c_{\alpha\bet}(x_0,\xi_0){y^\prime}^\alpha {\eta^\prime}^\beta.
\end{equation}
For the next condition set $\sigma^M_{x\xi}(P)=\sum_{j=0}^{M\mu/(\mu+1)}
\sum_{\alp/\mu+\beta=M-j(\mu+1)/\mu}c_{\alpha\beta}(x,\xi)(y^\prime)^{\alpha}D^\beta_{y^\prime}$.
\begin{equation}\label{TransCond2}
	\text{The kernel of } \sigma^M_{x\xi}(P)\text{ in } \cS(\R^n) \text{ is } \{0\}.
\end{equation}
Therefore in order to prove Theorem \ref{MainMetivier} it is enough to show the following statement.
\begin{Thm}\label{MetivierThm2}
	Let $P$ be of the form \eqref{canonicalform1} in a conic neighborhood $\cU$ of $(0,e_n)$.
	If $P$ satisfies \eqref{TransCond1} and \eqref{TransCond2} then there
	is a conic neighborhood $\crb$ of $(0,e_n)$ such 
	that 
	\begin{equation*}
		\WF_{[\fM]}u\cap\crb=\WF_{[\fM]} u\cap\crb\qquad \fa u\in\E^\prime(\Omega)
	\end{equation*}
	for any [semiregular] weight matrix $\fM$.
\end{Thm}
In order to prove Theorem \ref{MetivierThm2} we need to continue to transform the operator $P$.
We fix the following notation
\begin{equation*}
	A_j=\frac{\partial}{\partial x_j},\qquad A_{-j}=x_j\left(\frac{\partial}{\partial x_n}\right)^{\tfrac{1}{\mu}},\qquad j=1,\dotsc,\nu.
\end{equation*}
For $I=(j_1,\dotsc,j_k)\in\Set*{\pm 1,\dotsc,\pm\nu}^{k}$ we set $A_I=A_{j_1}\dotsc A_{j_k}$.
Moreover we denote $\abs{I_+}=\#\{j_k>0\}$, $\abs{I_-}=\#\{j_k<0\}$ 
and $\langle I\rangle=\abs{I_+}+(1/\mu)\abs{I_-}$.
Thus if $P$ is of the form \eqref{canonicalform1} then we can write
\begin{equation*}
	P(x,D_x)=\sum_{\langle I\rangle=M}c_{I}(x,D_x)A_I
\end{equation*}
where $c_I$ are analytic pseudodifferential operators of order $d-M$ in a neighborhood of $(x_0,\xi_0)$.
If we multiply $P$ with an elliptic factor and, if necessary, take powers of $P$
then we can assume that $d=M>\nu$.

As the next step we add new variables $x^{\dprime}=(x_{-\nu},\dotsc,x_{-1})$ and
and denote $\tilde{x}=(x^{\dprime},x)$ and $\tilde{\xi}=(\xi^\dprime,\xi)$.
Let $\phi(x^\dprime)\in\CC^\infty_0(\R^\nu)$ be such that $\phi(x^{\dprime})=1$ near the origin.
We extend a distribution $u\in\D^\prime(\R^n)$ to a distribution $\tilde{u}\in\D^\prime(\R^{n+\nu})$
by setting
\begin{equation*}
	\tilde{u}(\tilde{x})=\phi(x^\dprime)u(x).
\end{equation*}
The operators $A_j$ and $A_{-j}$ are extended by
\begin{equation*}
	\tilde{A}_j=\frac{\partial}{\partial x_j},\qquad \tilde{A}_{-j}=\left(\frac{\partial}{\partial x_{-j}}+x_j\frac{\partial}{\partial x_n}\right)
	\left(\frac{\partial}{\partial x_n}\right)^{-(\mu-1)/\mu}.
\end{equation*}
If we consider $c_I(x,\xi)$ as a symbol in some conic neighborhood of $(\tilde{x}_0,\tilde{\xi_0})$, $\tilde{x}_0=(0,x_0)$, $\tilde{\xi}_0=(0,\xi_0)$,
which is independent of $(x^\dprime,\xi^\dprime)$ then we can extend the operator $P$ by
\begin{equation*}
	\tilde{P}(\tilde{x},D_{\tilde{x}})=\sum_{\langle I\rangle=M}\tilde{c}_I(\tilde{x},D_{\tilde{x}})\tilde{A}_I
\end{equation*}
We observe that there are a neighborhood $U$ of $x_0$ 
and a conic neighborhood $\tilde{\cU}$ 
of $(\tilde{x},\tilde{\xi})$
such that for all $\tilde{u}\in\E^\prime(U)$ we have that
\begin{equation}\label{ExtendWF}
	\WF_{[\fM]} \widetilde{Pu}\cap\cU=\WF_{[\fM]}\tilde{P}\tilde{u}\cap\cU.
\end{equation}
Finally we consider the coordinate change $\tilde{x}\mapsto \tilde{y}=(y^\dprime,y)$ given by
\begin{align*}
	y^\dprime&=(y_{-\nu},\dotsc,y_{-1})=(x_{-\nu},\dotsc,x_{-1})\\
	y&=(y_1,\dotsc,y_{n-1},y_n)=\left(x_1,\dotsc,x_{n-1},x_n-\frac{1}{2}\sum_{j=1}^\nu x_jx_{-j}\right)
\end{align*}
In the $\tilde{y}$-coordinates the operator $\tilde{P}$ transfroms to
\begin{equation*}
	Q\left(\tilde{y},D_{\tilde{y}}\right)=\sum_{\langle I\rangle =M}d_I(\tilde{y},D_{\tilde{y}})
	X_I
\end{equation*}
where the $d_I$ are operators of order $0$ and
\begin{align*}
	X_j&=\frac{\partial}{\partial y_j}-\frac{1}{2}y_{-j}\frac{\partial}{\partial y_{n}}\\
	X_{-j}&=\left(\frac{\partial}{\partial_{-y_j}}+\frac{1}{2}y_j\frac{\partial}{\partial y_n}\right)
	\left(\frac{\partial}{\partial y_n}\right)^{-(\mu-1)/(\mu)}.
\end{align*}
In this setting Okaji \cite{MR0807494} proved the following theorem:
\begin{Thm}[{\cite[Theorem 3]{MR0807494}}]\label{OkajiThm}
	Let $N=n+\nu$ and $\cU\subseteq T^\ast\R^N\setminus\{0\}$ be a conic neighborhood 
	of $(x_0,\xi_0)$, $x_0=0$ and $\xi_0=(0,\dotsc,0,1)$. Assume that the operator $P$
	is defined in $\cU$ and satisfies the following conditions: 
	\begin{itemize}
		\item $P(x,D_x)=\sum_{\langle I\rangle=M}c_I(x,D_x)X_I$ where $c_I$ are classical
		pseudodifferential operators of order $0$ in $\cU$, 
		$X_j=\frac{\partial}{\partial x_j}-\frac{1}{2}x_{j+\nu}\frac{\partial}{\partial x_N}$
		and $X_{-j}=(\frac{\partial}{\partial x_{j+\nu}}+\frac{1}{2}x_j\frac{\partial}{\partial x_N})\left(\frac{\partial}{\partial x_N}\right)^{-(\mu-1)/\mu}$ for $j=1,\dotsc,\nu$ 
		and $M\geq \nu+  1$.
		\item For any $\zeta\in\R^{2\nu}\setminus\{0\}$ we have 
		$\sum_{\langle I\rangle=M}c_{I,0}(x_0,\xi_0)\zeta^I\neq 0$ where $c_{I,0}$ denotes
		the principal symbol of the operator $c_I$.
		\item Setting $\mathit{P}_{x\xi}(y,D_y)=\sum_{\langle I\rangle=M}\tilde{X}_I(y,D_y)$;
		$\tilde{X}_j=\frac{\partial}{\partial y_j}-\frac{1}{2}y_{j+\nu}$ and 
		$\tilde{X}_{-j}=\frac{\partial}{\partial y_{j+\nu}}+\frac{1}{2}y_j$ then
		$\mathit{P}_{x_0\xi_0}(y,D_y)$ is injective on $\cS(\R^N)$.
	\end{itemize}
	Then there are a neighborhood $U$ of $x_0$, a conic neighborhood $\cU$ of $(x_0,\xi_0)$
	and an operator $A\in\op(S^{-M/\mu}_{1/(\mu+1),1/(\mu+1)})(U)$ such that 
	for all $\varphi\in\CC^\infty_0(U)$ and every $u\in\E^\prime(U)$ we have that
	\begin{equation}\label{Okaji1}
		\cU\cap\WF_a(A\varphi u)=\emptyset.
	\end{equation}
\end{Thm}
In Theorem \ref{OkajiThm} $\op(S^{m}_{\rho,\delta}(\Omega))$ denotes
a class of analytic pseudodifferential operators of type $(\rho,\delta)$, $0<\rho\leq 1$ and 
$0\leq\delta<1$, introduced by
Metivier \cite{MR0597752}.
A function $a\in \CC^\infty(U\times U\times\R^N)$, $U\subseteq\R^N$ open, 
is in $a\in S^m_{\rho,\delta}$ if there are $C>0$ and $R>0$ such that
\begin{equation*}
	\abs*{\partial^{\alpha}_{x,y}\partial^\beta_\xi a(x,y,\xi)}\leq C^{\alp+\bet+1}
	\left(1+\xit\right)^m\left(\alp+\alp^{1-\delta}\xit^\delta\right)^\alp\left(\frac{\bet}{\xit}\right)^{\rho\bet}
\end{equation*}
for all $\alpha\in\Z_+^{2N},\beta\in\Z_+^N$, $x,y\in U$ and $\xi\in\R^N$ with $R\bet\leq\xit$.
Then the kernel of the associated operator $A=\OPA a$ is
\begin{equation*}
	K_A=(2\pi)^{-N}\int e^{i(x-y)\xi}a(x,y,\xi)\,d\xi.
\end{equation*}
According to \cite[Lemma 3.3]{MR0597752}, cf.~also \cite[pp.~33-34]{MR0597752} we have that
\begin{equation*}
	\WF_a(K_A)\subseteq\Set*{(x,y;\xi,\eta)\given x=y,\; \xi=-\eta}.
\end{equation*}
Thence by Theorem \ref{WFkernel}
\begin{equation}\label{Delta-Local}
	\WF_{[\fM]}Au\subseteq\WF_{[\fM]}u
\end{equation}
for all $u\in\E^\prime(U)$ and every [semiregular] weight matrix $\fM$.
We conclude also from \eqref{Okaji1} that
\begin{equation}\label{Okaji2}
	\cU\cap\WF_{[\fM]} A\varphi u=\emptyset.
\end{equation}
Therefore combining \eqref{ExtendWF}, \eqref{Delta-Local} and \eqref{Okaji2}
proves Theorem \ref{MetivierThm2}.
\appendix
\section{Differential operators with constant coefficients}
\begin{Thm}\label{ConstantThm}
	Let $P=P(D)$ be an operator with constant coefficients
	in $\R^n$ 
	and $\fM$ be a [semiregular] weight matrix.
	Then the following statements are equivalent:
	\begin{enumerate}
		\item The operator $P$ is $[\fM]$-microhypoelliptic
		in any open set $\Omega\subseteq\R^n$, i.e.~$\WF_{[\fM]}Pu=\WF_{[\fM]}u$ for all $u\in\D^\prime(\Omega)$.
		\item The operator $P$ is $[\fM]$-hypoelliptic in 
		any open set $\Omega\subseteq\R^n$.
		\item For any open set $\Omega\subseteq\R^n$ we have that
		\begin{equation*}
			\Set*{u\in\D^\prime(\Omega)\given Pu=0}\subseteq\DC{\fM}{\Omega}.
		\end{equation*}
		\item Every fundamental solution $E$ of $P$ satisfies
		\begin{equation*}
			\singsupp_{[\fM]} E=\{0\}.
		\end{equation*}
		\item There exists a fundamental solution $E$ of $P$
		such that 
		\begin{equation*}
			\singsupp_{[\fM]} E=\{0\}.
		\end{equation*}
	\end{enumerate}
\end{Thm}
\begin{proof}
	The implications $(1)\Rightarrow(2)\Rightarrow(3)\Rightarrow(4)\Rightarrow (5)$
	are trivial.
	Suppose now that $(5)$ holds and let $v\in\E^\prime(\R^n)$.
	Corollary \ref{ConvWF} implies that 
	\begin{equation}\label{ConstantMicro1}
\WF_{[\fM]} v=	\WF_{[\fM]} \bigl(E\ast(Pv)\bigr)=\WF_{[\fM]}Pv.
	\end{equation}
	Let  $u\in\D^\prime(\Omega)$ be a distribution in some open set $\Omega\subseteq\R^n$.
	We continue by  choosing  an arbitrary point $p$ and a neighborhood 
	$U\Subset\Omega$. If $\chi\in\CC^\infty_0(\Omega)$ is a test function
	with $\chi\vert_U=1$ then $(P(\chi u))\vert_U=(Pu)\vert_U$.
	Hence \eqref{ConstantMicro1} implies for any $\xi\in\R^n\!\setminus\!\{0\}$ that
	\begin{equation*}
		(p,\xi)\in\WF_{[\fM]} u\Longleftrightarrow (p,\xi)\in\WF_{[\fM]}Pu.
	\end{equation*} 
	It follows that $P$ is $[\fM]$-microhypoelliptic in
	$\Omega$.
\end{proof}

We will say that $P$ is $[\fM]$-hypoelliptic if
one of the equivalent conditions of Theorem \ref{ConstantThm}
is satisfied for some [semiregular] weight matrix $\fM$.
It is easy to see that any $[\fM]$-hypoelliptic operator with constant
coefficients is smooth hypoelliptic,
cf.~\cite[Theorem 2.1]{Treves2006}.
More generally, the characterization in
Theorem \ref{ConstantThm} shows that $[\fM]$-hypoellipticity implies
hypoellipticity for any larger ultradifferentiable class:
\begin{Cor}
	Let $P$ be a differential operator with constant coefficients.
	\begin{itemize}
		\item If $P$ is $[\fM]$-hypoelliptic
		for some [semiregular] weight matrix 
		then $P$ is $[\fN]$-hypoelliptic for all
		[semiregular] weight matrices $\fN$ such that
		$\bM[\preceq]\fN$.
		\item If $P$ is $\{\fM\}$-hypoelliptic for some
		R-semiregular weight matrix then $P$
		is $(\fN)$-hypoelliptic for any B-semiregular weight matrices $\fN$ with $\bM\{\lhd)\bN$.
		\item If $P$ is $(\fM)$-hypoelliptic for some B-semiregular weight matrix $\fM$ then
		$P$ is $\{\fN\}$-hypoelliptic for all R-semiregular weight matrices $\fN$ such that
	$\fM\{\preceq\}\fN$.
	\end{itemize}
\end{Cor}
The following Theorem is a slight generalization of \cite[Theorem 11.4.7]{Hoermander2005}.
We may write $\bM\preceq\bN$ for two weight sequences if there is a constant $C>0$ such that
$M_k\leq C^{k+1}N_k$.
\begin{Thm}\label{ConstantThm2}
	Let $y\in\R^n$, $P$ be a hypoelliptic differential operator with constant coefficients, $\fM$ be a countable collection
	of weight sequences,
	$\Omega\subseteq\R^n$ an open set and $x_0\in\Omega$.
	We assume that for any $u\in\CC^\infty(\Omega)$
	with $Pu=0$ there are some $h>0$ and $\bM\in\fM$ such that
	\begin{equation*}
		\abs*{\langle y,D\rangle^k u(x_0)}\leq h^k M_k,\qquad k\in\Z_+^n.
	\end{equation*}
	Then there are $\bM^0\in\fM$, $s\geq 1$ and $C>0$ such that 
	$\bG^s\preceq \bM^0$ and 
	\begin{equation}\label{ZeroEst}
		\abs{\langle y,\zeta\rangle}\leq C\bigl(1+\abs*{\imag \zeta}\bigr)^s
	\end{equation}
	for all $\zeta\in\C^n$ satisfying $P(\zeta)=0$.
\end{Thm}
\begin{proof}
	Let 
	\begin{equation*}
		\mathcal{N}=\Set*{u\in L^2_{loc}(\Omega)\given Pu=0}.
	\end{equation*}
	Since $P$ is hypoelliptic we have that $\mathcal{N}\subseteq\CC^\infty(\Omega)$ and
	using \cite[Theorem 4.4.2]{MR1996773} we conclude that
	the topologies on $\mathcal{N}$ induced by $L^2_{loc}(\Omega)$ and by $\CC^\infty(\Omega)$
	coincide.
	For $r\in\N$ and $\bM\in\fM$
	we set
	\begin{equation*}
		F_{r,\bM}=\Set*{u\in\mathcal{N}\given \abs*{\langle y,D\rangle^k u(x_0)}\leq r^k M_k\;\fa k\in\N}
	\end{equation*}
	 which is a closed subset of the Frechet space $\mathcal{N}$.
	 By assumption $\bigcup_{r\in\N,\bM\in\fM} F_{r,\bM}=\mathcal{N}$.
	 Baire's Theorem implies that there have to be some $r_0\in\N$
	 and a weight sequence $\bM^0\in\fM$ such that
	 $F_{r_0,\bM^0}$ includes an interior point.
	 Since $F_{r_0,\bM^0}$ is convex and symmetric, the origin
	 has to be an interior point. That means
	 there have to be some $\delta>0$ and a compact subset
	 $K$ of $\Omega$ such that $F_{r_0,\bM^0}$ includes
	 all functions $u\in\mathcal{N}$ with $\norm{u}{L^2(K)}\leq \delta$.
	 Then 
	 \begin{equation}\label{ConstantEst}
	 	\abs*{\langle y,D\rangle^k u(x_0)}\leq r_0^{k} M^0_{k}
	 	\delta^{-1}\norm{u}{L^2(K)},\qquad \fa u\in\mathcal{N}.
	\end{equation}
	Observe that \eqref{ConstantEst} is homogeneous in $u$
	and obviously true for all $u\in\mathcal{N}$ with
	$\norm{u}{L^2(K)}=\delta$ and therefore for all $u\in\mathcal{N}$.
	
	Now we can follow the end of the proof of
	\cite[Theorem 11.4.7]{Hoermander2005} verbatim
	to conclude that there is some $s\geq 1$
	and $C>0$ such that
	\begin{equation*}
		\abs*{\langle y,\zeta\rangle}\leq C\left(1+\abs*{\imag \zeta}\right)^s
	\end{equation*}
	for all $\zeta\in\C^n$ with $P(\zeta)=0$
	and $\bG^s\preceq\bM^0$.
	%Now let $\zeta\in\C$ with $P(\zeta)=0$ and set
%	$u(x)=e^{i\langle x,\zeta\rangle}$ in \eqref{ConstantEst}.
%	If $A$ is a bound of the function $x\mapsto\abs{x-x_0}$ on $K$ then \eqref{ConstantEst}  implies that there is a constant $C>0$ such that
%	\begin{equation*}
	%	\abs*{\langle y,\zeta\rangle}^k\leq %Cr_0^kM_k^0 e^{A\abs{\imag\zeta}}
%	\end{equation*}
	%%for all $k\in\N$.
	%	Let $\mu(\tau)=\sup\Set*{\abs*{\langle y,\zeta\rangle}\given P(\zeta)=0,\; \abs{\imag \zeta}\leq \tau}$.
	%	As in the proof Theorem
	%	gives that there are $s,C\geq 1$ such that
	%	$\mu(t)= 2C\tau^s(1+o(1))$
\end{proof}

\begin{Cor}\label{ConstantCorFin}
	Let $P$ be a differential operator with constant coefficients
	and $\fM$ be a weight matrix.
Then the following holds:
\begin{enumerate}
	\item If $\bM$ is R-semiregular and $P$ is $\{\fM\}$-hypoelliptic then there exists $s\geq 1$ 
	such that $P$ is $\{\bG^s\}$-hypoelliptic and
	$\{\bG^s\}\{\preceq\}\fM$.
	\item If $\bM$ is B-semiregular and $P$ is $(\fM)$-hypoelliptic then there exists $s\geq 1$
	such that $P$ is $\{\bG^s\}$-hypoelliptic and 
	$\{\bG^s\}\{\lhd)\fM$.
\end{enumerate}
\end{Cor}

\begin{proof} 
	If $\fM$ is $R$-semiregular and $P$ is $\{\fM\}$-hypoelliptic
	then $\Set*{u\in L^2_{loc}(\Omega)\given Pu=0}\subseteq\Rou{\fM}{\Omega}$ by Theorem \ref{ConstantThm}.
	Hence we can apply Theorem \ref{ConstantThm2} and conclude
	for each $j=1,\dotsc,n$ there are
	$\bM^j\in\fM$, $s_j\geq 1$ and $C_j$ such that
	\begin{equation*}
		\abs*{\langle e_j,\zeta\rangle}\leq C_j
		\left(1+\abs*{\imag \zeta}\right)^{s_j}
	\end{equation*}
for all zeros $\zeta$ of the polynomial $P(\zeta)$
and $\bG^{s_j}\preceq\bM^{j}$.

Now let $s=\max s_j$ and $\bM^0$ the largest
weight sequence of  the collection $\Set{\bM^j\given j=1,\dotsc,n}$.\footnote{Note that by definition $\fM$ is totally ordered in the sense of the pointwise order.}
It follows first that $\bG^s\preceq\bM^0$
which by definition means that $\{\bG^s\}\{\preceq\}\fM$.

Moreover
\begin{equation*}
	\abs{\zeta}\leq \sum_{j=1}^n\abs*{\zeta_j}\leq
	\sum_{j=1}^nC_j\left(1+\abs*{\imag \zeta}\right)^{s_j}\leq C\left(1+\abs*{\imag \zeta}\right)^s
\end{equation*}
where $C:=\sum C_j$. Hence \cite[Theorem 2.2.2]{MR1249275} implies that
$P$ is $\{\bG^s\}$-hypoelliptic.
	
	In the Beurling case we need to return to the proof of Theorem \ref{ConstantThm2}. By assumption $P$ is $(\fM)$-hypoelliptic
	and therefore
	\begin{equation}\label{Zero}
		\mathcal{N}=\Set*{u\in L^2_{loc}(\Omega)\given Pu=0}\subseteq\Beu{\fM}{\Omega}=\bigcap_{\bM\in\fM}\Beu{\bM}{\Omega}.
	\end{equation}
In particular, $P$ is also hypoelliptic 
and therefore for each $j\in\{1,\dotsc,n\}$ there is some 
$s_j\geq 1$ such that \eqref{ZeroEst} holds for $y=e_j$ following the proof
of \cite[Theorem 11.4.7]{Hoermander2005} (see also \cite[Section 2.2]{MR1249275})
As in the proof of the Roumieu case we conclude that
$P$ is $\{\bG^{s_0}\}$-hypoelliptic for $s_0=\max s_j$.
Let $j_0$ be such that \eqref{ZeroEst} holds for $y=e_{j_0}$
only for $s_0$.
From \eqref{Zero} we infer that for all $u\in\mathcal{N}$,
$h> 0$ and every $\bM\in\fM$ there is some $C>0$ such that
\begin{equation*}
	\abs*{D_{j_0}^ku(x_0)}\leq Ch^kM_k,\qquad k\in\N.
\end{equation*}
Setting $u=e^{i\langle x,\zeta\rangle}$ for some $\zeta$ with
$P(\zeta)=0$ we obtain that
for all $h>0$ and $\bM\in\fM$ there is some $C>0$ such that
\begin{equation}\label{Zeta}
	\abs*{\zeta_{j_0}}^k\leq Ch^k M_k,\qquad k\in\Z_+. 
\end{equation}
As in the proof of \cite[Theorem 11.4.7]{Hoermander2005}
we set $\mu(\tau)=\max\Set{\abs{\zeta_{j_0}}\given \zeta\in\C^n:\; P(\zeta)=0\,\&\, \abs{\imag (\zeta)\leq\tau}}$
and recall that by the definition of $s_0$ the proof of
\cite[Theorem 11.4.7]{Hoermander2005} gives that there is some
$b>0$ such that $\mu(\tau)\geq b\tau^{s_0}$ for large $\tau$.
Therefore \eqref{Zeta} implies that
\begin{equation*}
	\tau^{s_0 k}\leq C\left(\frac{h}{b}\right)^kM_k
\end{equation*}
for large $\tau$. Putting $\tau=k$ and using Stirling's estimates we see that
for each $\bM\in\fM$ and $h>0$ there is some $C>0$
such that
\begin{equation*}
	(k!)^{s_0}\leq Ch^k M_k
\end{equation*}
for large $k$. Thence $\bG^{s_0}\{\lhd)\fM$.
\end{proof}
If we consider the Gevrey matrix $\fG=\Set{\bG^s\given s>1}$
then it follows directly that a differential operator with 
constant coefficients is hypoelliptic
if and only if the operator is $\{\fG\}$-hypoelliptic.
This assertion extends trivially to all classes which are larger than
$\Rou{\fG}{\Omega}$:
\begin{Cor}
	Let $P$ be a differential operator with constant coefficients and $\fM$ be a weight matrix.
	Then the following holds:
	\begin{itemize}
		\item If $\fM$ is R-semiregular and $\fG\{\preceq\}\fM$
		then $P$ is hypoelliptic if and only if
		$P$ is $\{\fM\}$-hypoelliptic.
		\item If $\fM$ is B-semiregular and $\fG\{\lhd)\fM$
		then $P$ is hypoelliptic if and only if
		$P$ is $(\fM)$-hypoelliptic.
	\end{itemize}
\end{Cor}
We conclude by considering the quasianalytic case.
We recall that,
if $\fM$ is a weight matrix then 
$\Rou{\fM}{\Omega}$ is quasianalytic if and only if 
any $\bM\in\fM$ is a quasianalytic weight sequence, i.e.~$\bM$
does not satifies \eqref{NQ}.
In turn, $\Beu{\fM}{\Omega}$ is quasianalytic, if and only if
there is some $\bM\in\fM$ which is quasianalytic.

If $\fM$ is an $R$-semiregular weight matrix with $\Rou{\fM}{\Omega}$ being quasianalytic 
and $P$ is $\{\fM\}$-hypoelliptic then necessarily $s=1$ in Theorem \ref{ConstantCorFin}(1), i.e.~$P$ is analytic hypoelliptic.
On the other hand if $\fM$ is $B$-semiregular, $\Beu{\fM}{\Omega}$ is quasianalytic and 
$P$ is $(\fM)$-hypoelliptic then
$P$ is also analytic hypoelliptic by
Corollary \ref{ConstantCorFin}(2).
Since any differential operator $P$ with constant coefficients is analytic hypoelliptic exactly when $P$ is elliptic we
conclude that the following statement holds, which
generalizes the classical Theorem of Petrowsky \cite{Petrowsky1939} regarding analytic hypoellipticity.
\begin{Cor}
	Let $P$ be a differential operator with constant coefficients
	and $\fM$ be a [semiregular] weight matrix such that
	$\DC{\fM}{\R^n}$ is quasianalytic.
	Then $P$ is $[\fM]$-hypoelliptic if and only if $P$ is elliptic.
\end{Cor}

\section{A semiregular quasianalytic class\\ transversal to all semiregular Denjoy-Carleman classes}
The following theorem is a generalization of the construction in \cite[Subsection 7.4]{MR4002151}. 
In this appendix $\E^{\{\ast\}}$, $\ast=\bM,\fM,\omega$, refers
to the subsheaf of smooth functions generated by $\bM$, $\fM$ or $\omega$ (in the Roumieu sense).
\begin{Thm}\label{TransversalThm}
	Let $\fM=\Set*{M^\lambda\given \lambda\in \N}$ be a ordered countable family of sequences $M^\lambda=(M_k^{(\lambda)})_{k\in\Z_+}$ of positive numbers such that
	\begin{enumerate}
		\item $\sup_{k\in\N}\sqrt[k]{M_k^{(\lambda)}}=+\infty$ for all $\lambda\in\N$.
		\item $M_k^{{(\lambda_1)}}\leq M_k^{(\lambda_2)}$ for 
		all $k\in\Z_+$ and all $\lambda_1\leq\lambda_2$.
	\end{enumerate}
	Then there exists a weight function $\omega$ in the sense of Braun-Meise-Taylor with the following properties:
	\begin{itemize}
	\item	$\omega$ is semiregular
	(cf.~Definition \ref{Def:Semiregular}(2)).
	\item $\E^{\{\omega\}}$ is quasianalytic.
	\item $\E^{\{\omega\}}\nsubseteq\E^{\{\fM\}}$.
\end{itemize}
\end{Thm}

\begin{proof}
	We start by defining several auxiliary sequences:
	First we select two sequences $(a_\lambda)_{\lambda\in\N}$, $(b_\lambda)_{\lambda\in\N}$ iteratively.
	First we choose $a_1\in\N$ to be the smallest number such
	that $({M^{(1)}_{a_1}})^{2/a_1}\geq 1$ and set
	\begin{equation*}
		\left(M_{a_1}^{(1)}\right)^{2/a_1}
	\end{equation*}
	 which is possible by (1).
	Moreover, we iteratively select $a_k\in\N$ in such a way that
	\begin{equation*}
		b_k:=\left(M_{a_k}^{(k)}\right)^{2/a_k}
	\end{equation*}
	satisfies the conditions $b_{k+1}\geq b_k$ and $b_k\geq k^2$
	and 
	\begin{equation*}
		a_{k+1}\geq e^{b_k}a_k+1
	\end{equation*} 
	for all $k\in\N$. This is again possible by (1).
	It follows that $k\mapsto a_k$ is strictly increasing and
	therefore $\lim a_k=+\infty$.
	
	We set
	\begin{equation*}
		q_k:=1, \quad 0\leq k <a_1,\quad q_k:=b_{\lambda}^k,\quad
		a_\lambda\leq k<a_{\lambda+1},\quad \lambda\in\N.
	\end{equation*}
	We may set $A_\lambda:=\Set{k\in\N\given a_\lambda\leq k<a_{\lambda+1}}$ and $\rho_k=q_k/q_{k-1}$.
	%\newline
	
	\textbf{Claim:} The sequence $\rho_k$ is non-decreasing.\newline
	In fact, if $k-1,k\in A_\lambda$ for some $\lambda$.
	then $\rho_k=b_\lambda$ and therefore $\rho_k$ is constant
	on the set $B_\lambda:=\Set{k\in\N\given k-1,k\in A_\lambda}$. If $k-1\in A_{\lambda-1}$ and $k\in A_{\lambda}$ then $\rho_{k}=b^k_\lambda/b^{k-1}_{\lambda-1}$.
	Then $\rho_{k+1}=b_\lambda\geq \rho_k$ since $b_{\lambda-1}\geq b_\lambda$ and analogously
	$\rho_{k-1}= b_{\lambda-1}\leq b_\lambda^k/b_\lambda^{k-1}$.
	Note that $\lim a_\lambda=\infty$.
	
	Hence $\mathbf{q}=(q_k)_k$ is a weight sequence
	and the same is true for $\mathbf{Q}=(k!q_k)_k$.
	In fact $\lim_{k\rightarrow\infty} \sqrt[k]{q_k}=\infty$,
	by construction, i.e. condition \eqref{DC-AnalyticIncl} holds for $\mathbf{Q}$.
	Thus $\An(\Omega)\subseteq\DC{\mathbf{Q}}{\Omega}$.
	
	%\newline
	\textbf{Claim:} $\E^{\{\mathbf{Q}\}}\nsubseteq\E^{\{\fM\}}$.
	\newline
	Let $\lambda\in\N$ then by definition we have for any $\nu\geq \lambda$ that
	\begin{equation*}
		\left(\frac{Q_{a_\nu}}{M_{a_k}^{(\lambda)}}\right)^{1/a_\nu}
		\geq \left( \frac{q_{a_\nu}}{M_{a_\nu}^{(\nu)}}\right)^{1/a_\nu}
		=\frac{b_k}{\bigl(M_{a_\nu}^{a_\nu}\bigr)^{1/a_\nu}}
		=\left(M_{a_\nu}^{(\nu)}\right)^{1/a_\nu}\geq \nu.
	\end{equation*}
	It follows that
	\begin{equation}
		\sup_{k\in\N}\left(\frac{Q_k}{M_k^{\lambda}}\right)=\infty
	\end{equation}
	 and therefore $\E^{\{\mathbf{Q}\}}\nsubseteq\E^{\{\bM^{(\lambda)}\}}$ for any
	$\lambda\in\N$ by \cite[Prop.~2.12(i)]{MR3285413}.
	Using also \cite[Prop.~4.6(i)]{MR3285413}
	we see that $\E^{\{\mathbf{Q}\}}\nsubseteq\E^{\{\fM\}}$.
	Note here that in order to apply these statements we
	need only that $\mathbf{Q}$ is a weight sequence
	but not necessarily that any $\bM^{(\lambda)}$ satisfies
	\eqref{LogConvex}.
	
	\textbf{Claim:} The weight sequence $\mathbf{Q}$ 
	is quasianalytic.
	
	Recall the asymptotic formula for the harmonic series
	$\sum_{k=1}^p=\log p+\gamma+\eps_p$ where $\gamma$ is
	the Euler constant and $\eps_p\rightarrow 0$ if $p\rightarrow\infty$.
	Thus
	\begin{equation*}
		\sum_{\ell\in A_\lambda}\frac{1}{\lambda}
		=\sum_{\ell=1}^{a_{\lambda+1}-1}\frac{1}{\ell}
		-\sum_{\ell=1}^{a_\lambda-1}\frac{1}{\ell}
		=\log (a_{\lambda+1}-1)-\log (a_\lambda-1)+\eps_{a_{\lambda+1}-1}-\eps_{a_\lambda-1},
	\end{equation*}
	and therefore the following estimate holds for any $\lambda$
	by construction:
	\begin{equation*}
		\sum_{\ell\in A_\lambda}\frac{1}{(Q_\ell)^{1/\ell}}
		\geq \sum_{\ell\in A_\lambda}\frac{\ell}{(q_\ell)^{1/\ell}}
		=\frac{1}{b_\lambda}\sum_{\ell\in A_\lambda}\frac{1}{\ell}.
	\end{equation*}
	On the other hand by the choice of $a_{\lambda+1}$
	we see that
	\begin{equation*}
		\frac{1}{b_\lambda}\log \left(\frac{a_{\lambda+1}-1}{a_\lambda-1}\right)\geq 1
	\end{equation*}
	for any $\lambda\geq 2$.
	Thence
	\begin{equation}\label{nq2}\tag{$\mathbf{nq}^\prime$}
		\sum_{k=0}^{\infty}\frac{1}{(Q_k)^{1/k}}=\infty
	\end{equation}
	and thus $\mathbf{Q}$ is quasianalytic, because
	\eqref{nq2} holds exactly when the series in \eqref{NQ}
	is divergent, cf.~e.g.~\cite{MR1996773} or \cite[Theorem Theorem 19.11]{Rudin1987}.
	
	If we consider now
	the weight function $\omega_{\mathbf{Q}}$ associated
	to the weight sequence $\mathbf{Q}$, given by \eqref{omegaM},
	then it is well known that $\omega_{\mathbf{Q}}$ 
	satisfies all conditions in Definition \ref{Def:BMT}
	except \eqref{omega1}, cf.~\cite[Lemma 2.3]{Schindl2021}
	and the literature cited there.
	However, since the sequence $\mathbf{Q}$ is strongly log-convex,
	i.e.~the sequence $\mathbf{q}$ satisfies \eqref{LogConvex},
	the associated weight function $\omega_{\mathbf{Q}}$
	is subadditive
	by \cite[Theorem 4.5]{Schindl2021}.
	This implies that $\omega_{\mathbf{Q}}$ satisfies \eqref{omega1}, i.e.~$\omega_{\mathbf{Q}}$ is a
	weight function in the sense of Braun-Meise-Taylor.
	Moreover, $\omega_{\mathbf{Q}}$ is semiregular:
	According to \cite[Lemma 2.3]{Schindl2021}
	the property $\lim_{k\rightarrow\infty} (q_k)^{1/k}=\infty$
	implies that $\omega_{\mathbf{Q}}(t)=o(t)$ for $t\rightarrow\infty$.
	
	Since $\E^{\{\mathbf{Q}\}}\subseteq\E^{\{\omega_{\mathbf{Q}}\}}$,
	cf.~\cite{Schindl2021}, it follows
	that $\E^{\{\omega_{\mathbf{Q}}\}}\nsubseteq\E^{\{\fM\}}$.
	Finally, the quasianalyticity of $\mathbf{Q}$ implies
	that 
	\begin{equation*}
		\int_{1}^{\infty}\frac{\omega_{\mathbf{Q}}(t)}{1+t^2}\,dt=\infty,
	\end{equation*}
	see \cite[Lemma 4.1]{Komatsu1973},
	which characterizes the quasianalyticity of
	 $\E^{\{\omega_{\mathbf{Q}}\}}$, cf.~\cite[Section 5]{MR3285413}.
\end{proof}
\begin{Cor}\label{TransversalCor}
	There exists a weight function $\omega$ in the sense of Braun-Meise-Taylor such that $\E^{\{\omega\}}$ is quasianalytic and
	\begin{equation*}
		\E^{\{\omega\}}\nsubseteq\E^{\{\bM\}}
	\end{equation*}
	for any semiregular weight sequence $\bM$.
\end{Cor}
\begin{proof}
	If $\bM$ is a semiregular weight sequence then
	\eqref{DC-DerivClosed} gives that there is a constant $\rho>0$ such that $M_k\leq \rho^{k^2+k}$.
	Thence $\bM\preceq\bN^\rho$, where $\bN^\rho$ is the 
	weight sequence given by $N^\rho_k=\rho^{k^2}$.
	If $\fN:=\Set{\bN^\rho\given \rho>0}$ then this implies
	that $\E^{\{\bM\}}\subseteq\E^{\{\fN\}}$ for any semiregular weight sequence $\bM$.
	 However, it is easy to see that $\E^{\{\fN\}}=\E^{\{\widehat{\fN}\}}$ where
	 $\widehat{\fN}=\Set{\bN^\rho\given \rho\in\N}$
	 and that $\fM=\widehat{\fN}$ satisfies the conditions
	 in Theorem \ref{TransversalThm}.
	 By Theorem \ref{TransversalThm} there exists therefore
	 a semiregular weight function $\omega$ such that
	 $\E^{\{\omega\}}$ is quasianalytic and
	 $\E^{\{\omega\}}\nsubseteq\E^{\{\fN\}}$.
\end{proof}
\begin{Rem}
	It is worth pointing out that by construction
	the weight function $\omega$ in Theorem \ref{TransversalThm}
	resp.~Corollary \ref{TransversalCor} is not only semiregular
	but also subadditive.
	For example, that means that $\E^{\{\omega\}}$ can be 
	characterized by almost-analytic extensions according to
	\cite[Theorem 1.5]{MR4002151} and
	that it satisfies the assumptions of \cite[Theorem 4.8]{Albanese2010} and \cite[Theorem 7.7]{MR4002151}, cf.~\cite[Theorem 4.8]{MR4002151}.
\end{Rem}
\bibliographystyle{plainurl}
\bibliography{hypoelliptic}

\begin{thebibliography}{10}

\bibitem{Albanese2010}
A.~A. Albanese, D.~Jornet, and A.~Oliaro.
\newblock Quasianalytic wave front sets for solutions of linear partial
  differential operators.
\newblock {\em Integral Equations Operator Theory}, 66(2):153--181, 2010.
\newblock \href {https://doi.org/10.1007/s00020-010-1742-6}
  {\path{doi:10.1007/s00020-010-1742-6}}.

\bibitem{Albanese2012}
A.~A. Albanese, D.~Jornet, and A.~Oliaro.
\newblock Wave front sets for ultradistribution solutions of linear partial
  differential operators with coefficients in non-quasianalytic classes.
\newblock {\em Math. Nachr.}, 285(4):411--425, 2012.
\newblock \href {https://doi.org/10.1002/mana.201010039}
  {\path{doi:10.1002/mana.201010039}}.

\bibitem{Albano2013}
Paolo Albano and Antonio Bove.
\newblock Wave front set of solutions to sums of squares of vector fields.
\newblock {\em Mem. Amer. Math. Soc.}, 221(1039):xii+73, 2013.
\newblock \href {https://doi.org/10.1090/S0065-9266-2012-00663-0}
  {\path{doi:10.1090/S0065-9266-2012-00663-0}}.

\bibitem{Bonet2007}
José Bonet, Reinhold Meise, and Sergej~N. Melikhov.
\newblock {A comparison of two different ways to define classes of
  ultradifferentiable functions}.
\newblock {\em Bulletin of the Belgian Mathematical Society - Simon Stevin},
  14(3):425--444, 2007.
\newblock \href {https://doi.org/10.36045/bbms/1190994204}
  {\path{doi:10.36045/bbms/1190994204}}.

\bibitem{BoutetdeMonvel1976}
Louis Boutet~de Monvel, Alain Grigis, and Bernard Helffer.
\newblock Parametrixes d'op\'erateurs pseudo-diff\'erentiels \`a{}
  caract\'eristiques multiples.
\newblock In {\em Journ\'ees: \'Equations aux {D}\'eriv\'ees {P}artielles de
  {R}ennes (1975)}, volume No. 34--35 of {\em Ast\'erisque}, pages 93--121.
  Soc. Math. France, Paris, 1976.

\bibitem{Bove2018}
Antonio Bove and Gregorio Chinni.
\newblock Analytic and {G}evrey hypoellipticity for perturbed sums of squares
  operators.
\newblock {\em Ann. Mat. Pura Appl. (4)}, 197(4):1201--1214, 2018.
\newblock \href {https://doi.org/10.1007/s10231-017-0720-x}
  {\path{doi:10.1007/s10231-017-0720-x}}.

\bibitem{Bove2024}
Antonio Bove and Marco Mughetti.
\newblock Minimal {G}evrey regularity for {H}\"ormander operators.
\newblock {\em Int. Math. Res. Not. IMRN}, (4):2790--2832, 2024.
\newblock \href {https://doi.org/10.1093/imrn/rnad055}
  {\path{doi:10.1093/imrn/rnad055}}.

\bibitem{Bove2013}
Antonio Bove, Marco Mughetti, and David~S. Tartakoff.
\newblock Hypoellipticity and nonhypoellipticity for sums of squares of complex
  vector fields.
\newblock {\em Anal. PDE}, 6(2):371--445, 2013.
\newblock \href {https://doi.org/10.2140/apde.2013.6.371}
  {\path{doi:10.2140/apde.2013.6.371}}.

\bibitem{MR1052587}
R.~W. Braun, R.~Meise, and B.~A. Taylor.
\newblock Ultradifferentiable functions and {F}ourier analysis.
\newblock {\em Results Math.}, 17(3-4):206--237, 1990.
\newblock URL: \url{http://dx.doi.org/10.1007/BF03322459}, \href
  {https://doi.org/10.1007/BF03322459} {\path{doi:10.1007/BF03322459}}.

\bibitem{BraunRodrigues2021}
Nicholas Braun~Rodrigues and Antonio~Victor da~Silva, Jr.
\newblock Approximate solutions of vector fields and an application to
  {D}enjoy-{C}arleman regularity of solutions of a nonlinear {PDE}.
\newblock {\em Math. Nachr.}, 294(8):1452--1471, 2021.
\newblock \href {https://doi.org/10.1002/mana.201800516}
  {\path{doi:10.1002/mana.201800516}}.

\bibitem{Chinni2023}
Gregorio Chinni.
\newblock On the regularity of the solutions and of analytic vectors for ``sums
  of squares''.
\newblock In {\em Bruno {P}ini {M}athematical {A}nalysis {S}eminar 2022},
  volume 13(1) of {\em Bruno Pini Math. Anal. Semin.}, pages 90--108. Univ.
  Bologna, Alma Mater Stud., Bologna, 2023.
\newblock \href {https://doi.org/10.6092/issn.2240-2829/16159}
  {\path{doi:10.6092/issn.2240-2829/16159}}.

\bibitem{Cordaro2024}
Paulo~D. Cordaro and Stefan F\"urd\"os.
\newblock The {M}etivier inequality and ultradifferentiable hypoellipticity.
\newblock {\em Math. Nachr.}, 297(7):2517--2531, 2024.
\newblock \href {https://doi.org/10.1002/mana.202300147}
  {\path{doi:10.1002/mana.202300147}}.

\bibitem{MR2126468}
Paulo~D. Cordaro and Nicholas Hanges.
\newblock Impact of lower order terms on a model {PDE} in two variables.
\newblock In {\em Geometric analysis of {PDE} and several complex variables},
  volume 368 of {\em Contemp. Math.}, pages 157--176. Amer. Math. Soc.,
  Providence, RI, 2005.
\newblock \href {https://doi.org/10.1090/conm/368/06777}
  {\path{doi:10.1090/conm/368/06777}}.

\bibitem{Derridj1973}
M.~Derridj and C.~Zuily.
\newblock Sur la r\'egularit\'e{} {G}evrey des op\'erateurs de {H}\"ormander.
\newblock {\em J. Math. Pures Appl. (9)}, 52:309--336, 1973.

\bibitem{MR226355}
Klaus Floret and Joseph Wloka.
\newblock {\em Einf\"uhrung in die {T}heorie der lokalkonvexen {R}\"aume},
  volume No. 56 of {\em Lecture Notes in Mathematics}.
\newblock Springer-Verlag, Berlin-New York, 1968.

\bibitem{MR4002151}
Stefan F\"{u}rd\"{o}s, David~Nicolas Nenning, Armin Rainer, and Gerhard
  Schindl.
\newblock Almost analytic extensions of ultradifferentiable functions with
  applications to microlocal analysis.
\newblock {\em J. Math. Anal. Appl.}, 481(1):123451, 51, 2020.
\newblock \href {https://doi.org/10.1016/j.jmaa.2019.123451}
  {\path{doi:10.1016/j.jmaa.2019.123451}}.

\bibitem{Fuerdoes2022}
Stefan F{\"u}rd{\"o}s and Gerhard Schindl.
\newblock The theorem of iterates for elliptic and non-elliptic operators.
\newblock {\em J. Funct. Anal.}, 283(5):74, 2022.
\newblock Id/No 109554.
\newblock \href {https://doi.org/10.1016/j.jfa.2022.109554}
  {\path{doi:10.1016/j.jfa.2022.109554}}.

\bibitem{Fuerdoes2024A}
Stefan F{\"u}rd{\"o}s and Gerhard Schindl.
\newblock Ellipticity and the problem of iterates in {Denjoy}-{Carleman}
  classes.
\newblock {\em Collect. Math.}, 2024.
\newblock \href {https://doi.org/10.1007/s13348-024-00455-7}
  {\path{doi:10.1007/s13348-024-00455-7}}.

\bibitem{Hoepfner2020}
G.~Hoepfner and R.~Medrado.
\newblock Microlocal regularity for {M}izohata type differential operators.
\newblock {\em J. Inst. Math. Jussieu}, 19(4):1185--1209, 2020.
\newblock \href {https://doi.org/10.1017/s1474748018000361}
  {\path{doi:10.1017/s1474748018000361}}.

\bibitem{Hoepfner2023}
G.~Hoepfner, A.~Raich, and P.~Rampazo.
\newblock Weighted {H}ardy spaces and ultradistributions.
\newblock {\em Houston J. Math.}, 49(1):195--229, 2023.

\bibitem{Hoermander1955}
Lars H\"ormander.
\newblock On the theory of general partial differential operators.
\newblock {\em Acta Math.}, 94:161--248, 1955.
\newblock \href {https://doi.org/10.1007/BF02392492}
  {\path{doi:10.1007/BF02392492}}.

\bibitem{MR0294849}
Lars H\"{o}rmander.
\newblock Uniqueness theorems and wave front sets for solutions of linear
  differential equations with analytic coefficients.
\newblock {\em Comm. Pure Appl. Math.}, 24:671--704, 1971.
\newblock \href {https://doi.org/10.1002/cpa.3160240505}
  {\path{doi:10.1002/cpa.3160240505}}.

\bibitem{MR1996773}
Lars H\"{o}rmander.
\newblock {\em The analysis of linear partial differential operators. {I}}.
\newblock Classics in Mathematics. Springer-Verlag, Berlin, 2003.
\newblock Distribution theory and Fourier analysis, Reprint of the second
  (1990) edition [Springer, Berlin; MR1065993 (91m:35001a)].
\newblock \href {https://doi.org/10.1007/978-3-642-61497-2}
  {\path{doi:10.1007/978-3-642-61497-2}}.

\bibitem{Hoermander2005}
Lars H\"{o}rmander.
\newblock {\em The analysis of linear partial differential operators. {II}}.
\newblock Classics in Mathematics. Springer-Verlag, Berlin, 2005.
\newblock Differential operators with constant coefficients, Reprint of the
  1983 original.
\newblock \href {https://doi.org/10.1007/b138375} {\path{doi:10.1007/b138375}}.

\bibitem{Komatsu1973}
Hikosaburo Komatsu.
\newblock Ultradistributions. {I}. {S}tructure theorems and a characterization.
\newblock {\em J. Fac. Sci. Univ. Tokyo Sect. IA Math.}, 20:25--105, 1973.

\bibitem{MR551623}
Gottfried K\"othe.
\newblock {\em Topological vector spaces. {II}}, volume 237 of {\em Grundlehren
  der Mathematischen Wissenschaften}.
\newblock Springer-Verlag, New York-Berlin, 1979.
\newblock \href {https://doi.org/10.1007/978-1-4684-9409-9}
  {\path{doi:10.1007/978-1-4684-9409-9}}.

\bibitem{Liess1999}
Otto Liess.
\newblock Carleman regularization in the {${C}^\infty$}-category.
\newblock volume~45, pages 213--240. 1999.
\newblock Workshop on Partial Differential Equations (Ferrara, 1999).
\newblock \href {https://doi.org/10.1007/BF02826096}
  {\path{doi:10.1007/BF02826096}}.

\bibitem{Matsumoto1987}
Waichir\^o Matsumoto.
\newblock Theory of pseudodifferential operators of ultradifferentiable class.
\newblock {\em J. Math. Kyoto Univ.}, 27(3):453--500, 1987.
\newblock \href {https://doi.org/10.1215/kjm/1250520659}
  {\path{doi:10.1215/kjm/1250520659}}.

\bibitem{Metivier1980}
Guy M\'etivier.
\newblock Une classe d'op\'erateurs non hypoelliptiques analytiques.
\newblock {\em Indiana Univ. Math. J.}, 29(6):823--860, 1980.
\newblock \href {https://doi.org/10.1512/iumj.1980.29.29059}
  {\path{doi:10.1512/iumj.1980.29.29059}}.

\bibitem{MR0597752}
Guy M\'{e}tivier.
\newblock Analytic hypoellipticity for operators with multiple characteristics.
\newblock {\em Comm. Partial Differential Equations}, 6(1):1--90, 1981.
\newblock \href {https://doi.org/10.1080/03605308108820170}
  {\path{doi:10.1080/03605308108820170}}.

\bibitem{MR0807494}
Takashi \={O}kaji.
\newblock Analytic hypoellipticity for operators with symplectic
  characteristics.
\newblock {\em J. Math. Kyoto Univ.}, 25(3):489--514, 1985.
\newblock \href {https://doi.org/10.1215/kjm/1250521068}
  {\path{doi:10.1215/kjm/1250521068}}.

\bibitem{Okaji1988}
Takashi \=Okaji.
\newblock Gevrey-hypoelliptic operators which are not
  {$C^\infty$}-hypoelliptic.
\newblock {\em J. Math. Kyoto Univ.}, 28(2):311--322, 1988.
\newblock \href {https://doi.org/10.1215/kjm/1250520484}
  {\path{doi:10.1215/kjm/1250520484}}.

\bibitem{Petrowsky1939}
I.~G. Petrowsky.
\newblock Sur l'analyticite des solutions des syst{\`e}mes d'{\'e}quations
  diff{\'e}rentielles.
\newblock {\em Rec. Math. Moscou, n. Ser.}, 5:3--68, 1939.

\bibitem{MR3285413}
Armin Rainer and Gerhard Schindl.
\newblock Composition in ultradifferentiable classes.
\newblock {\em Studia Math.}, 224(2):97--131, 2014.
\newblock \href {https://doi.org/10.4064/sm224-2-1}
  {\path{doi:10.4064/sm224-2-1}}.

\bibitem{MR1249275}
Luigi Rodino.
\newblock {\em Linear partial differential operators in {G}evrey spaces}.
\newblock World Scientific Publishing Co., Inc., River Edge, NJ, 1993.
\newblock \href {https://doi.org/10.1142/9789814360036}
  {\path{doi:10.1142/9789814360036}}.

\bibitem{Rudin1987}
Walter Rudin.
\newblock {\em Real and complex analysis}.
\newblock McGraw-Hill Book Co., New York, third edition, 1987.

\bibitem{Schindl2021}
Gerhard Schindl.
\newblock On subadditivity-like conditions for associated weight functions.
\newblock {\em Bull. Belg. Math. Soc. Simon Stevin}, 28(3):399--427, 2021.
\newblock \href {https://doi.org/10.36045/j.bbms.210127}
  {\path{doi:10.36045/j.bbms.210127}}.

\bibitem{MR0296509}
Fran\c{c}ois Tr\`eves.
\newblock Analytic-hypoelliptic partial differential equations of principal
  type.
\newblock {\em Comm. Pure Appl. Math.}, 24:537--570, 1971.
\newblock \href {https://doi.org/10.1002/cpa.3160240407}
  {\path{doi:10.1002/cpa.3160240407}}.

\bibitem{MR0597144}
Fran\c{c}ois Tr\`eves.
\newblock {\em Introduction to pseudodifferential and {F}ourier integral
  operators. {V}ol. 1}.
\newblock University Series in Mathematics. Plenum Press, New York-London,
  1980.
\newblock Pseudodifferential operators.

\bibitem{MR0597145}
Fran\c{c}ois Tr\`eves.
\newblock {\em Introduction to pseudodifferential and {F}ourier integral
  operators. {V}ol. 2}.
\newblock University Series in Mathematics. Plenum Press, New York-London,
  1980.
\newblock Fourier integral operators.

\bibitem{Treves2006}
Fran\c{c}ois Tr\`eves.
\newblock {\em Basic linear partial differential equations}.
\newblock Dover Publications, Inc., Mineola, NY, 2006.
\newblock Reprint of the 1975 original.

\bibitem{MR4436039}
Fran\c{c}ois Tr\`eves.
\newblock {\em Analytic partial differential equations}, volume 359 of {\em
  Grundlehren der mathematischen Wissenschaften [Fundamental Principles of
  Mathematical Sciences]}.
\newblock Springer, Cham, 2022.
\newblock \href {https://doi.org/10.1007/978-3-030-94055-3}
  {\path{doi:10.1007/978-3-030-94055-3}}.

\end{thebibliography}
\end{document}